\documentclass[11pt]{article} 
\title{\bf Geometric Anosov flows of dimension five with smooth distributions} 
\author{Yong Fang} 
\date{\it Laboratoire de Math\'ematique d'Orsay, U.M.R. 
8628 du C.N.R.S, Universit\'e Paris-Sud, France\\ (e-mail: fangyong1@yahoo.fr)} 

\usepackage{amsmath,amsthm,amsfonts,amssymb,amscd,amstext} 
\chardef\bslash=`\\

\theoremstyle{definition}

\theoremstyle{remark}

\begin{document} 
\maketitle 
\renewcommand{\sectionmark}[1]{} 
{\bf Abstract}--{\it We classify the five dimensional $C^\infty$ Anosov flows 
which have $C^\infty$-Anosov splitting and preserve a smooth 
pseudo-Riemannian metric . Up to a special time change and 
finite covers, such a flow is $C^\infty$ flow equivalent either to 
the suspension of a symplectic hyperbolic automorphism of 
$\mathbb{T}^{4}$, or to the geodesic flow on a three dimensional 
hyperbolic manifold.}\\\\
1. Introduction\\
2. Preliminaries

    2.1. Some generalities

    2.2. Proof of Theorem 1\\
3. Homogeneity in dimension $5$ 

    3.1. Remarks about rank $0$ and $4$

    3.2. Homogeneity in rank $2$\\
4. The case of $2$ positive Lyapunov exponents

    4.1. Preparations

    4.2. dim$\mathfrak{h}'=1$

    4.3. dim$\mathfrak{h}'=2$ \\
5. The case of 1 positive exponent and $d\lambda\wedge\omega\not\equiv 0$

    5.1. Preparations

    5.2. dim$\mathfrak{h}' =2$

    5.3. dim$\mathfrak{h}'=1$\\
6. The case of $1$ positive exponent and $d\lambda\wedge\omega\equiv 0$

    6.1. Preparations
 
    6.2. dim$\mathfrak{h}'=1$

    6.3. dim${\mathfrak{h}}'=2$ \\
7. Appendix\\\\
{\bf 1.  Introduction}

Let $M$ be a $C^{\infty}$-closed manifold. A $C^{\infty}$-flow, $\phi_{t}$, generated 
by the non-singular vector field $X$ is called an Anosov flow if there exists a $\phi_t$-invariant 
splitting of the tangent bundle $$ TM= \mathbb{R}X\oplus E^{+}\oplus E^{-},$$
a Riemannian metric on $M$ and two positive numbers $a$ and $b$, such that
$$\forall \  u^{\pm}\in E^{\pm}, \  \forall \  t \geq 0,  \  \parallel D\phi_{\mp t}(u^{\pm})\parallel\leq
a e^{-bt}\parallel u^{\pm}\parallel ,$$
where $E^{-}$ and $E^{+}$ are called   
the strong stable and strong unstable distributions of the flow.

In general , $E^{-}$ and $E^{+}$ are only continuous. 
If they are both $C^{\infty}$ subbundles of $TM$, then the Anosov flow is said to have smooth 
distributions. This case is rather rare, see for example [$\mathbf{K}$], 
[$\mathbf{FK}$], and [$\mathbf{BFL 2}$].  Although the smoothness of these 
two distributions is dynamically so strong a condition, it is still quite weak
geometrically. So to arrive at a classification result, one has to suppose 
in addition the existence of a smooth invariant geometric structure.
For example, in {[\bf BFL2]}, the existence of an invariant 
contact form is assumed. 

If an Anosov flow preserves a $C^\infty$ pseudo-Riemannian metric, then by definition, 
this flow is called {\it geometric}. In this paper, we consider the {\it geometric} Anosov flows 
with smooth distributions.

The classical examples of such flows are the suspensions of symplectic hyperbolic
infranilautomorphisms and the geodesic flows 
on locally symmetric spaces of rank one. There exist also lots of 
non-classical algebraic models (see {\bf [To]}), which makes a possible classification of such flows 
quite interesting. In this 
paper, we obtain the classification in dimension five.

In general, given an Anosov flow with $C^{\infty}$ distributions $\phi_{t}$, one gets 
a smooth $1$-form $\lambda$, such that
$$ \lambda( E^\pm) =0, \  \lambda (X) =1.$$
It is called the {\it canonical 1-form} of the flow,which is easily seen
to be $\phi_{t}$-invariant.\\\\
{\bf Definition.}  ${\rm rank}(\phi_{t}):= 2({\rm max} \{ k\geq 0 \mid \wedge ^{k}d\lambda\not \equiv 0 \})$.

$\ $ 

We call this even number the {\it rank} of $\phi_{t}$. Here $\wedge ^{k}d\lambda $ denotes the exterior 
$k$-th power of $d\lambda $, and by convention, $\wedge ^{0}d\lambda $ $:=1.$ Note that rank($\phi_{t}$)
is just the rank of the $2$-form $d\lambda $ (see {\bf [Lich]}). If $\phi_t$ is topologically transitive and 
its rank is $2k$, then $\wedge^k d\lambda$ vanishes nowhere on an open-dense subset of $M$.  

For $\forall\  a\in \mathbb{R}$, denote by $[a]$ the biggest integer, which is smaller than $a$. 
If the dimension of $M$ is $m$, 
then the degree of $\wedge ^{[\frac{m}{2}]+1}d\lambda$ will be bigger than $m$. So we have
$$ {\rm rank}(\phi_{t})\leq 2 [\frac{m}{2}].$$
In Section $2$, we characterize 
the classical homogeneous models above by their ranks. More precisely, we prove

$\ $

{\bf Theorem 1.} {\it Let $M$ be a $C^{\infty}$ closed manifold of dimension $m$ 
and $\phi_{t}$ be a geometric Anosov flow with $C^{\infty}$ distributions on $M$, we have\\
(i) if rank$(\phi_{t}) =0$, then up to a constant change of time scale, $\phi_{t}$ is 
$C^{\infty}$ flow equivalent to the suspension of a hyperbolic infranilautomorphism ;\\
(ii) if rank$(\phi_{t}) = 2[\frac{m}{2}]$, then up to finite covers, $\phi_{t}$ is 
$C^{\infty}$ flow equivalent to a canonical perturbation of the 
geodesic flow on a locally symmetric Riemannian manifold of strictly negative curvature.} 

$\ $

A {\it canonical perturbation} of a smooth flow with generator $X$ is (by definition)
the flow of the field $\frac{X}{1+\alpha(X)}$, where $\alpha$ is a $C^{\infty}$
closed $1$-form such that $1+\alpha(X) >0.$ It should be mentioned that 
Theorem $1.$ is just a more or less direct reformulation of the results of
{\bf [BFL2]}, {\bf [BL]} and {\bf [Plan]}.

Although there exist algebraic models of {\it geometric} Anosov flows 
with rank between $0$ and $2 [\frac{m}{2}]$, none of them is of dimension five. In fact, we prove the following

$\ $

{\bf Theorem 2.} {\it Let $M$ be a closed manifold of dimension five 
and $\phi_{t}$ be a geometric Anosov flow  with $C^{\infty}$ distributions 
on $M$, then\\
(i) either, up to a constant change of time scale and finite covers, $\phi_{t}$ is $C^{\infty}$ 
flow equivalent to the suspension of a symplectic hyperbolic 
automorphism of $\mathbb{T}^{4}$;\\
(ii)  or, up to finite covers, $\phi_{t}$ is $C^{\infty}$ flow equivalent to 
 a canonical perturbation of the geodesic flow on a three dimensional 
Riemannian manifold of constant negative curvature.}

$\ $

In the appendix, two lemmas are proved, which are used in the proof of 
Theorem $2$. Lemma $A$ is about the completeness of a linear connection and 
Lemma $B$ is about the time change of an Anosov flow with $C^\infty$ distributions.

If $M$ admits a {\it geometric} Anosov flow, then the dimension of 
$M$ must be odd (see Section $2$). In dimension three, an Anosov flow 
with $C^{\infty}$ distributions is {\it geometric} iff it preserves a  volume form (see {\bf[HK]}). 
Such flows are classified by \'E. Ghys (see {\bf [Gh]}). Here Theorem 2. gives 
a classification for the case of dimension five. We should mention that 
such five dimensional flows are also studied in {\bf [FK]} with the 
purpose to understand the contact case. 

Beginning with
dimension seven, we can find many algebraic models of {\it geometric} 
Anosov flows, which are neither contact nor suspensions (see {\bf [To]}). 
The situation will then become much more complex and 
a classification is still out of reach at the moment. Indeed, our proof of Theorem $2.$ 
is quite specific to the case of dimension five.\\\\ 
{\bf 2.  Preliminaries.}\\
{\bf 2.1.  Some generalities.}  

Let $\phi_{t}$ be an Anosov flow with $C^\infty$ distributions on a $C^\infty$ closed manifold $M$. 
Denote by $X$ the generator of this flow. For each $C^\infty$ $2$-form $\omega$ on $M$, denote by 
Ker$\omega$ the kernel of $\omega$, i.e. Ker$\omega :=\{ y\in TM\mid i_y\omega =0\}.$ Let us prove 
at first the following

$\ $
 
{\bf Lemma 2.1.1.} {\it Under the notations above, $\phi_t$ is geometric, iff it preserves a $C^\infty$ 
$2$-form with $\mathbb{R}X$ as kernel.}

{\it Proof.} Suppose that $\phi_t$ is {\it geometric}. Denote by $g$ a $C^\infty$ $\phi_t$-invariant 
pseudo-Riemannian metric. Then by the Anosov property of $\phi_t$, we get 
$$ g(X, E^\pm )=0,\  g(E^\pm, E^\pm)=0.$$
Let $J$ be the section of $T^\ast M\otimes TM$, such that 
$$ J(X) =0,\  J(u^\pm)= \pm u^\pm,\  \forall\  u^\pm\in E^\pm.$$
Then $g(J\cdot, \cdot)$ is easily seen to be a $C^\infty$ $\phi_t$-invariant $2$-form, denoted by 
$\omega$. Since $g$ is non-degenerate, then so is $\omega\mid_{E^+\oplus E^-}$. 
Again by the Anosov property, we get $i_X\omega =0.$ So the kernel of 
$\omega$ is $\mathbb{R}X$. 

Suppose that $\phi_t$ preserves a $C^\infty$ $2$-form $\Theta$, such that Ker$\Theta =\mathbb{R}X.$ Then 
there exists a unique $\phi_t$-invariant symmetric $(0, 2)$-tensor $g$, such that 
$$ g(X, X) =1,\  g(X, u^\pm) =0,$$
$$ g(u^+, u^-) = g(u^-, u^+)=\Theta (u^+, u^-),$$
$$ g(u^\pm, v^\pm) =0, \  \forall \  u^\pm,\  v^\pm\in E^\pm.$$
Since Ker$\Theta =\mathbb{R}X$, then $g$ is non-degenerate. So $g$ is a pseudo-Riemannian metric. 
Thus $\phi_t$ is {\it geometric}.  $\square$

$\ $

We deduce that the following Anosov flows 
with $C^\infty$ distributions are {\it geometric} :\\
(i). Contact Anosov flows with $C^\infty$ distributions. \\
(ii). Suspensions of symplectic hyperbolic infranilautomorphisms.\\
(iii). Three dimensional volume preserving Anosov flows with $C^\infty$ 
distributions. (see {\bf[HK]}).

In {\bf [To]}, P. Tomter constructed explicitly a seven dimensional Anosov flow, which is indeed {\it geometric}. 
By generalizing his ideas, we can then construct many non-usual algebraic models 
of {\it geometric} Anosov flows. The following lemma gives another way to construct such flows

$\ $

{\bf Lemma 2.1.2.} {\it Under the notations above, if $\phi_t$ is geometric, then for 
each $C^\infty$ $1$-form $\beta$, such that $\mathcal{L}_X d\beta =0$ and 
$\beta(X) >0,$ the flow of $\frac{X}{\beta(X)}$ is also a geometric Anosov flow with $C^\infty$ distributions.}

{\it Proof.} Denoted by $\phi^{\beta}_t$ the flow of $\frac{X}{\beta(X)}.$ Then by Lemma $B$ 
proved in the appendix, $\phi^{\beta}_t$ is also an Anosov 
flow with $C^\infty$ distributions.

Since $\phi_t$ is {\it geometric}, then by Lemma $2.1.1$, it preserves a $C^\infty$ $2$-form $\omega$, 
such that Ker$\omega =\mathbb{R}X.$ In particular, we have $i_X\omega =0.$ Then 
$$ i_Xd\omega =\mathcal{L}_X\omega -di_X\omega =0.$$
Thus
$$ \mathcal{L}_{X_\beta}\omega = i_{X_\beta}d\omega +
di_{X_\beta}\omega =0.$$
So $\phi^{\beta}_t$ preserves also $\omega$ and 
Ker$\omega =\mathbb{R}{X_\beta}.$ Then by Lemma $2.1.1$, $\phi^{\beta}_t$ is also {\it geometric}.  $\square$

$\ $

Let $\phi_t$ be as above and {\it geometric}. Since $\phi_t$ preserves a $C^\infty$ $2$-form $\omega$, 
such that Ker$\omega=\mathbb{R}X,$ then $\omega\mid_{E^+\oplus E^-}$ is non-degenerate. 
By the Anosov property of $\phi_t$, we get $\omega(E^\pm, E^\pm) =0$. So $E^+$ and $E^-$ are 
both Lagrangian subspaces of $\omega\mid_{E^+\oplus E^-}$. We deduce that $E^+$ and $E^-$ have the same 
dimension, denoted by $n$. So the dimension of $M$ is odd. 

It is easily seen that $\lambda\wedge(\wedge^n\omega)$ is a $\phi_t$-invariant volume form. So $\phi_t$ is 
topologically transitive (see {\bf [HaK]}). Denote by $\nu$ the probability defined by this volume form. Then 
by the Multiplicative Ergodic
Theorem of Oseledec, there exists a $\nu$-conull 
$\phi_t$-invariant subset $\Lambda$ of $M$ and a 
decomposition of $TM\mid_\Lambda$ into $\phi_{t}$-invariant 
measurable subbundles, 
$$  TM\mid_{\Lambda} = \oplus_{ 0\leq i\leq k} L_{i} ,$$
such that for $\forall$ $u_{i}\in L_{i}$,  
$$\lim_{ t\to \pm\infty} t^{-1}\log\parallel D\phi_t(u_{i})\parallel =\chi_{i},$$
where $L_{i}$ is called a Lyapunov subbundle and $\chi_{i}$ its Lyapunov 
exponent. $L_{i}$ is also denoted by $L_{\chi_{i}}$. 

The following lemma is due to Feres and Katok 
(see {\bf [FK]}).

$\ $

{\bf Lemma 2.1.3.} {\it Under the notations above, if 
$\tau$ is a $C^\infty$ $\phi_{t}$-invariant tensor 
field of type 
$(0, r)$ and $\sum_{ 1\leq l\leq r} \chi_{i_l} \not = 0$, then 
$\tau ( L_{i_1},\cdots, L_{i_r} ) =0.$}\\\\
{\bf 2.2.  Proof of theorem 1.}

Let $\phi_t$ be 
a {\it geometric} Anosov flow with $C^\infty$ distributions and 
suppose that $E^+$ is of dimension $n$. 
Then by the previous subsection, we have 
$m =2n+1,$ where $m$ is the dimension of $M$.

If rank$(\phi_t) =2[\frac{m}{2}]\  (= 2n),$ then $\wedge^n d\lambda\not
\equiv 0.$ So the $\phi_t$-invariant $C^\infty$ $m$-form 
$\lambda\wedge(\wedge^n d\lambda)$ is not 
identically zero. Since $\phi_t$  is topologically transitive, 
then $\exists\  c\not =0$, 
such that $\lambda\wedge(\wedge^n d\lambda) 
= c\cdot\lambda\wedge(\wedge^n \omega).$ 
We deduce that $ \lambda\wedge(\wedge^n d\lambda)$ vanishes nowhere, 
i.e. $\lambda$ is a contact form. 
Then by the classification of contact Anosov flows with $C^\infty$ distributions (see {\bf [BFL2]}), the 
case $(ii)$ of Theorem $1$. is true.

If  rank$(\phi_t) =0,$ then $d\lambda \equiv 0$. So $E^+\oplus E^-$ is 
integrable. By Theorem $3.1.$ of {\bf [Plan]}, $\phi_t$ admits 
a global section $\Sigma$ (a global section is by definition 
a connected closed submanifold of codimension $1$ which 
intersects each orbit transversally). Denote by $\tau$ the {\it first return time} function 
of $\Sigma$. Then the Poincar\'e map of $\Sigma$ 
is by definition $\psi :=\phi_{\tau(\cdot)}(\cdot).$ For the 
sake of completeness, we prove in detail  the following.

$\ $

{\bf Lemma 2.2.1.} {\it The previous Poincar\'e
 map $\psi$ is a $C^\infty$ Anosov diffeomorphism
with $C^\infty$ distributions, topologically 
transitive and preserving a $C^\infty$ linear connection.}

{\it Proof.} Recall that $E^+\oplus \mathbb{R} X$ and 
$E^-\oplus \mathbb{R}X$ are called the unstable and stable distributions of $\phi_t$.
They are both integrable (see {\bf [HaK]}). Denote by $\mathcal{F}^{+,0}$ and $\mathcal{F}^{-, 0}$ 
their corresponding foliations. Since $\Sigma$ is transversal to $X$, 
then $\mathcal{F}^{+,0}\cap 
\Sigma$ gives a $C^\infty$ foliation on $\Sigma$. Denote by $E_{\Sigma}^+$ 
its $C^\infty$ tangent distribution. Similarly we denote by $E_{\Sigma}^-$ the tangent distribution of 
$\mathcal{F}^{-, 0}\cap\Sigma$. 

Since $\mathcal{F}^{+, 0}$ is $\phi_t$-invariant, then the foliation 
$\mathcal{F}^{+,0}\cap\Sigma$ is $\psi$-invariant. We deduce that 
$E_\Sigma^+$ is $\psi$-invariant. Similarly $E_\Sigma^-$ is also $\psi$-invariant.

Fix a Riemannian metric on $M$. Since $E^+\mid_\Sigma$ 
and $E_\Sigma^+$ are both transversal to $\mathbb{R}X$ (along 
$\Sigma$), then we can project $E_\Sigma^+$ onto $E^+\mid_\Sigma$
with respect to $\mathbb{R}X$. Denote this projection by $P^+$. Since $\Sigma$ is compact, 
then we can find two positive constants $M_1$ and $M_2$, such that
$$M_{1}\parallel u\parallel\leq\parallel P^{+} u\parallel\leq M_{2}\parallel u\parallel,
\  \forall\  u\in E_{\Sigma}^{+}.$$
For $\forall\  x\in \Sigma$, take $u\in (E_{\Sigma}^+)_{x}$. Then $u$ splitts uniquely as 
$$u=P^{+}_{x}(u)+ a X_x,\  a\in\mathbb{R}.$$
We have
$$(D_{x}\psi)(u)= (D_{x}\tau(u)+ a) X_{\psi(x)}+(D_{x}\phi_{\tau(x)})(P_{x}^{+} u).$$
Thus
$$(D_{x}\psi)(u)=(P_{\psi(x)}^+)^{-1}[(D_{x}\phi_{\tau(x)})(P_{x}^{+} u)].$$
So for $\forall\  n\in \mathbb{N},$  
$$(D_{x}\psi^n)(u)=(P_{\psi^n(x)}^{+})^{-1}(D_{x}\phi_{\tau(x)+\cdots+\tau(\psi^{n-1}(x))})(P^{+}_{x} u).$$
We have a similar formula for $E_\Sigma ^-$. Now a simple estimation shows that
$\psi$ is an Anosov  diffeomorphism with $C^\infty$ distributions, $E_\Sigma ^+$ and $E_\Sigma ^-$. 

Since $\phi_t$ is {\it geometric}, then it preserves a $C^\infty$ $2$-form $\omega$ whose 
kernel is $\mathbb{R}X.$ Restrict $\omega$ to a $C^\infty$ 2-form $\omega_\Sigma$ on $\Sigma$.
Then using the fact that $i_X\omega =0$, $\omega_\Sigma$ is seen to be $\psi$-invariant.
Since $\omega_\Sigma$ is non-degenerate, then $\psi$ preserves a volume form. 
We deduce that $\psi$ is topological transitive.

Now a direct calculation shows the existence of a $C^\infty$ $\psi$-invariant 
connection $\nabla$ on $\Sigma$, such that
$$\nabla \omega_\Sigma =0, \  \nabla E^\pm_\Sigma\subseteq E^\pm_\Sigma,$$
$$\nabla_{Y^\pm}Y^\mp =P^\mp_\Sigma[Y^\pm, Y^\mp ], \  \forall\  Y^\pm \subseteq E^\pm_\Sigma.$$ $\square$

$\ $

By {\bf [BL]} and the previous lemma, $\psi$ is seen to be $C^\infty$-conjugate to a hyperbolic 
infranilautomorphism. Then by Corollary $3.5.$ of {\bf [Plan]}, 
the integral manifolds of $E^{+}\oplus E^{-}$ are compact. So we can take 
a leaf of $E^{+}\oplus E^{-}$ as $\Sigma$. With respect to this section, 
the {\it first return time} function is contant. Then Theorem $1.$ follows.\\\\
{\bf 3. Homogeneity in dimension 5.}\\
{\bf 3.1. Remarks about rank 0 and 4.}

Now we begin to prove Theorem $2$. 
Suppose that $\phi_t$ satisfies the conditions in Theorem $2.$ Denote by $X$ the generator of 
$\phi_t$ and by $\nu$ its invariant volume form. By Lemma $2.1.1$, $\phi_t$ preserves a 
$C^\infty$ $2$-form $\omega$, such that Ker$\omega =\mathbb{R}X$, i.e. $\omega\mid_{E^+\oplus E^-}$ 
is non-degenerate. Thus by Lemma $2.1.3$, if $a$ is a Lyapunov exponent of $\phi_t$ with respect to 
$\nu$, then so is $-a$. Since $M$ is of dimension five, 
then there exist only two possibilities for the Lyapunov exponents of $\phi_t$,
$$(i). -a<0 <a,$$
$$(ii). -a< -b <0 <b<a.$$

$\ $

{\bf Lemma 3.1.1.} {\it Under the notations above, we have $d\omega \equiv 0.$}

{\it Proof.} Since $\omega$ is $\phi_t$-invariant, then
$$\mathcal{L}_X\omega =0,\  i_X\omega =0.$$
So
$$ i_Xd\omega =\mathcal{L}_X\omega -di_X\omega =0,$$
i.e.
$$ d\omega(X,\cdot,\cdot) \equiv 0.$$

If $\phi_t$ has only one positive Lyapunov exponent, i.e. 
the case $(i)$ above is true, then by Lemma $2.1.3$, $d\omega\equiv 0$.

If the case $(ii)$ above is verified, then the Lyapunov 
subbundles are all of dimension one. Again by Lemma $2.1.3$, 
$d\omega \equiv 0$.  $\square$

$\ $

The rank of $\phi_t$ can only be $0,$ $2,$ or $4$. If rank$(\phi_t) =4$, 
then by Theorem $1$, $\phi_t$ is finitely covered by a {\it canonical perturbation} of the geodesic flow 
on a three dimensional locally symmetric space of strictly negative curvature. But such a Riemannan space 
must have contant negative curvature. So Theorem $2.$ is true in this case.

If rank$(\phi_t) =0$, then by Theorem $1$, up to a constant change of time scale, $\phi_t$ is 
finitely covered by the suspension of a four dimensional hyperbolic nilautomorphism. But in dimension four, 
such a hyperbolic nilautomorphism must be $(\mathbb{T}^4,\  \bar A)$, where 
$\bar A$ is the induced application of an invertible hyperbolic matrix $A$ in $GL(4, \mathbb{Z}).$ By 
Lemma $3.1.1$, $\bar A$ is in addition symplectic. So Theorem $2.$ is true in this case. 

So to prove Theorem $2$, we need only prove the non-existence of the case 
of rank $2$. In the following, we suppose on the contrary that 
there exists a rank $2$ {\it geometric} Anosov flow $\phi_t$ with 
$C^\infty$ distributions on a closed five dimensional manifold $M$. 
In Subsection $3.2.$ below, this flow $\phi_t$ is proved to be homogeneous. Then in 
Sections $4,$ $5$ and $6$, all the possible homogeneous models are eliminated by
some dynamical and Lie theoretical arguments.\\\\
{\bf 3.2.  Homogeneity in rank 2.}

Denote by $\lambda$ the {\it canonical $1$-form} of $\phi_t$. Since 
rank$(\phi_t) =2$, then 
$$ d\lambda\not\equiv 0,\  d\lambda\wedge d\lambda \equiv 0.$$
Define $U:=\{ x\in M\mid (d\lambda)_x\not =0\}$. 
Since $\phi_t$ is topologically transitive and preserves 
$d\lambda$, then $U$ is a $\phi_t$-invariant open-dense subset of $M$. Denote by 
$\pi$ the projection of $TM$ onto $M$. We define 
$$E_1 :=\{ y\in E^+\oplus E^- \mid i_y d\lambda =0,\  \pi(y)\in U \}$$
and 
$$E_1^\pm := E_1\cap E^\pm.$$
Since $\phi_t$ preserves $d\lambda$, $E^+$ and $E^-$, then 
$E_1$, $E^+_1$ and $E^-_1$ are all $\phi_t$-invariant.

$\ $

{\bf Lemma 3.2.1.} {\it $E_1$ is a two dimensional $C^\infty$ subbundle of $TM\mid_U.$ 
$E_1^+$ and $E^-_1$ are both one dimensional $C^\infty$ subbundles of $TM\mid_U$. 
In addition, $E_1 =E^+_1\oplus E^-_1.$}

{\it Proof.} Since $d\lambda(X,\cdot) \equiv 0$, then we view $d\lambda$ as a section of $(E^+\oplus E^-)^\ast$. 
For $\forall \  x\in U,$ we have $(d\lambda)_x\not =0.$ So near $x$, we can find $C^\infty$ local sections of 
$E^+\oplus E^-$, $V_1$ and $V_2$, such that 
$$d\lambda(V_1,\  V_2)\equiv 1.$$
Denote by $V$ the $C^\infty$ local distribution spanned by $V_1$ and $V_2$ and denote by $V^\perp$ 
the orthogonal of $V$ with respect to $d\lambda\mid_{E^+\oplus E^-}$. 

Since $d\lambda\mid_V$ is non-degenerate, then 
$$V\cap V^\perp =\{0\}.$$
For $\forall\  u\in E^+\oplus E^-$, such that $\pi(u)$ near $x$, the following vector is contained in 
$V^\perp$,
$$P(u) := u-d\lambda(u, V_2(\pi(u)))\cdot V_1(\pi(u)) - d\lambda(V_1(\pi(u)), u)\cdot V_2(\pi(u)).$$
So we deduce that locally
$$E^+\oplus E^- =V\oplus V^\perp.$$
In addition, we see that the projection of 
$E^+\oplus E^-$ onto $V^\perp$ with respect to this 
direct sum decomposition is $C^\infty$. So $V^\perp$ must be also $C^\infty$. 

Since $d\lambda\mid_V$ is non-degenerate and $d\lambda\wedge d\lambda \equiv 0$, then 
$$d\lambda\mid_{V^\perp}\equiv 0.$$
Thus locally
$$E_1 =V^\perp.$$
In particular, $E_1$ is $C^\infty$ and two dimensional. Since $d\lambda(E^\pm, E^\pm) \equiv 0$, then 
for $\forall\  u\in E_1$, its projections to $E^+$ and $E^-$ are also contained in $E_1$. Thus 
$$E_1 =E^+_1\oplus E^-_1.$$
If for some $x$ in $U$, $(E^+_1)_x$ is of dimension two, then $(d\lambda)_x$ will be zero, which contradicts 
our assumption. Thus $E^+_1$ and $E^-_1$ are both of dimension one. In addition, they are 
evidently $C^\infty.$  $\square$
 
$\ $

{\bf Lemma 3.2.2.} {\it Under the notations above, the Lyapunov decomposition of 
$\phi_t$ is smooth.}

{\it Proof.} By definition, the Lyapunov decomposition of $\phi_t$ is called smooth, if 
there exists a $C^\infty$ decomposition of $TM$ and a $\phi_t$-invariant 
$\nu$-conull subset $\bar\Lambda$ of $M$, such that the Lyapunov decomposition 
is defined on $\bar\Lambda$ and coincides on $\bar\Lambda$ with this $C^\infty$ decomposition.

If $\phi_t$ has only one positive Lyapunov exponent, 
then its Lyapunov decomposition is just the restriction of that of Anosov onto a $\nu$-conull 
subset of $M$. Since $\phi_t$ has $C^\infty$ distributions, then the lemma is true in this case.

Suppose that $\phi_t$ has $2$ positive Lyapunov exponents $b < a$. Then 
there exists a $\nu$-conull subset $\Lambda$ of $M$, such that 
$$ TM\mid_\Lambda = L^+_1\oplus L^-_1\oplus L^+_2\oplus L^+_2
\oplus\mathbb{R}X,$$
where $L_1^\pm$ and $L^\pm_2$ are the Lyapunov subbundles with exponents $\pm b$ and $\pm a$ 
(see Subsections $2.1.$ and $3.1$). 

Since $U$ is a $\phi_t$-invariant open-dense subset and the flow is 
$\nu$-ergodic, then $U$ is $\nu$-conull. So $\nu (U\cap \Lambda)=1.$

Take $x\in U\cap \Lambda$ and $l^\pm_i\in (L^\pm_i)_x, \  i=1, 2.$
By Lemma $2.1.3$, we have 
$$d\lambda (l^+_1, l^-_2) =0,\  d\lambda(l^-_1, l^+_2) =0.$$ 
Since $(d\lambda)_x \not =0$, then we must have $d\lambda (l^+_1, l^-_1)\not =0$ or 
$d\lambda(l^+_2, l^-_2)\not =0.$

Suppose that $d\lambda (l^+_2, l^-_2)\not =0.$ Since 
$d\lambda\wedge d\lambda \equiv 0$, then we must have $d\lambda (l^+_1, l^-_1) =0$. 
So $l^+_1\in (E^+_1)_x$, i.e. $(L^+_1)_x =(E^+_1)_x .$ 
Similarly, we get $(L^-_1)_x = (E^-_1)_x.$

Since $\omega\mid_{E^+\oplus E^-}$ is non-degenerate 
and $\omega (l^+_1, l^-_2) =0$, 
then $\omega (l^+_1, l^-_1)\not =0$. We deduce that $(d\lambda \wedge \omega )_x \not =0.$ So 
$\lambda\wedge d\lambda\wedge\omega$ is not identically zero. 
Then by the topological transitivity of $\phi_t$, $\exists\  c\not =0$, such that 
$$\lambda\wedge d\lambda\wedge\omega =c\cdot\lambda\wedge\omega\wedge\omega.$$ 
So $\lambda\wedge d\lambda\wedge\omega$
is nowhere zero. We deduce that $d\lambda$ vanishes nowhere and $U =M$. In particular, 
$E_1$ and $E^\pm_1$ are all $C^\infty$ subbundles of $TM.$

So by the arguments above, for $\forall \  x\in \Lambda,$ $(E^\pm_1)_x = (L^\pm_1)_x$ or $(L^\pm_2)_x.$
Define
$$ \Lambda_i :=\{ y\in \Lambda\mid E^\pm_1(y) = L^\pm_i(y)\},\   i = 1, 2.$$
Then $\Lambda_1$ and $\Lambda_2$ are both mesurable and 
$\phi_t$-invariant. So one of them is $\nu$-conull. 
Suppose that $\nu(\Lambda_1) =1$. Then we have 
$E^\pm_1\mid_{\Lambda_1}= L^\pm_1\mid_{\Lambda_1}.$

By Lemma $2.1.3,$ we have on $\Lambda_1$, 
$$L^\pm_2 = [Ker (v\mapsto\omega(L^\mp_1, v))]\cap E^\pm.$$
Define two $\phi_t$-invariant $C^\infty$ subbundles of $TM$ as follows, 
$$ E^\pm_2 :=[ Ker(v\mapsto\omega(E^\mp_1, v))]\cap E^\pm.$$
Then we have $E^\pm_2\mid_{\Lambda_1} = L^\pm_2\mid_{\Lambda_1}.$ So the 
Lyapunov decomposition coincides on a conull set with a $C^\infty$ decomposition of $TM$. 

If $\nu(\Lambda_2) =1$, then similar argument works.  $\square$

$\ $

{\bf Remark 3.2.1.} If $\phi_t$ has two positive Lyapunov exponents, then by the proof of Lemma $3.2.2$, 
we have four $C^\infty$ line bundles on $M$, $E^\pm_1$ and $E^\pm_2$. We shall call
$$TM =\mathbb{R}X\oplus E^+_{1}\oplus E^-_{1} \oplus E^+_2\oplus E^-_2$$
the $C^\infty$ Lyapunov decomposition of $\phi_t$. The Lyapunov exponents of the 
corresponding Lyapunov subbundles of $E^\pm_{1, 2}$ are called respectively the Lyapunov exponents of 
$E^\pm_{1, 2}$. $E^\pm_i$ are also denoted by $E_{a^\pm_i}$, where $a^\pm_i$ are the Lyapunov exponents of 
$E^\pm_i$. If $a$ is not a Lyapunov exponent of $\phi_t$, then by convention, $E_a :=\{0\}.$    

If $\phi_t$ has only one positive Lyapunov exponent, then 
the $C^\infty$ Lyapunov decomposition of $\phi_t$ means $TM =\mathbb{R}X\oplus E^+\oplus E^-.$

$\ $
  
Now we can construct a $C^\infty$ connection $\nabla$, adapted to our situation.

If the flow has two positive Lyapunov
exponents, then there exists a unique $C^\infty$ connection $\nabla$ on $M$, such that
$$\nabla X =0,\  \nabla \omega =0,\  \nabla E^\pm_i\subseteq E^\pm_i,$$
$$\nabla_{Y^\pm_j}Y^\mp_i  = P^\mp_i [ Y^\pm_j, Y^\mp_i],\  \forall \  i, j\in \{1, 2\},$$
$$\nabla_X Y^\pm_i := [X, Y^\pm_i]\pm a_i Y^\pm_i, \  \forall\  Y^\pm_i\subseteq E^\pm_i,$$
where $a_i$ denotes the Lyapunov exponent of $E^+_i$ and $P^\pm_i$ 
represent the projections of $TM$ onto $E^\pm_i$. 

If $\phi_t$ has only one positive Lyapunov exponent $a$, 
then we get a similar $C^\infty$ connection $\nabla$, such that
$$\nabla X=0,\  \nabla\omega =0, \  \nabla E^\pm\subseteq E^\pm,$$
$$\nabla_{Y^\pm}Y^\mp =P^\mp [Y^\pm, Y^\mp],$$
$$ \nabla_X Y^\pm =[X, Y^\pm]\pm a Y^\pm,\  \forall\  Y^\pm\subseteq E^\pm,$$
where $P^\pm$ represent the projections of $TM$ onto $E^\pm$.
 
If a transformation of 
$M$ preserves $X$, $\omega$, and the $C^\infty$ Lyapunov decomposition, then it 
preserves also $\nabla$. In particular, $\nabla$ is $\phi_t$-invariant.

$\ $

{\bf Lemma 3.2.3.} {\it Under the notations above, if $K$ be a $C^\infty$ $\phi_t$-invariant tensor 
field of type $(1, l)$ on $M$, then $K( E_{a_1},\cdots , E_{a_l})\subseteq E_{a_1+\cdots +a_l}$, 
where $a_1,\cdots ,a_l$ are arbitrary Lyapunov exponents of $\phi_t$. In addition, we have 
$\nabla K=0$.}

{\it Proof.} By the same arguments as in Lemma $2.5.$ of {\bf [BFL1]}, we get for arbitrary 
Lyapunov exponents, $a_1,\cdots,a_l$, 
$$K( E_{a_1},\cdots , E_{a_l})\subseteq E_{a_1+\cdots +a_l}.$$   

Now let $Z_1,\cdots , Z_l$ be the sections of the smooth 
subbundles, $E_{a_1},\cdots ,E_{a_l}$. We have
$$(\nabla_X K)(Z_1, \cdots , Z_l) = \nabla_X(K(Z_1, \cdots, Z_l)) -
\sum_{1\leq i\leq l}K(Z_1,\cdots, \nabla_X Z_i,\cdots, Z_l)$$
$$= [X, K( Z_1,\cdots,Z_l)] + (\sum_{1\leq i\leq l}a_i)K( Z_1,\cdots, Z_l) - K([X, Z_1]+ a_1Z_1,\cdots)\cdots$$
$$=[X, K( Z_1,\cdots, Z_l)] - \sum_{1\leq i\leq l}K( Z_1,\cdots,[ X, Z_i],\cdots,  Z_l)$$
$$=(\mathcal{L}_X K)(Z_1,\cdots, Z_l) = 0. $$
So $\nabla_X K =0.$ Since $\nabla K$ is a $\phi_t$-invariant tensor of 
type $(1, l+1)$, then we have 
$$(\nabla_{E_{a_0}}K)(E_{a_1},\cdots,E_{a_l})\subseteq E_{a_0+\cdots +a_l}.$$
Since for $\forall\  a\in \mathbb{R},$ $\nabla E_a\subseteq E_a$, then 
$$(\nabla_{E_{a_0}}K)(E_{a_1},\cdots,E_{a_l})\subseteq E_{a_1+\cdots +a_l}.$$
So if $a_0\not =0$, we have $\nabla_{E_{a_0}}K =0.$ We deduce that $\nabla K =0.$  $\square$

$\ $

Denote by $T$ the torsion of $\nabla$ and by $R$ its curvature tensor. 
Then by the previous lemma, we have
$$\nabla T=0, \  \nabla R=0,\  T(E_{a_1}, E_{a_2})\subseteq E_{a_1 + a_2}\  ;$$
If $a_1 + a_2\not =0$, then 
$$ R(E_{a_1}, E_{a_2}) = 0.$$

Denote by $\widetilde M$ the universal cover of $M$ and by $\widetilde\nabla$ the 
lifted connection of $\nabla$. Then we have 

$\ $

{\bf Lemme 3.2.4.} {\it Under the 
notations above, the group of $\widetilde\nabla$-affine transformations of $\widetilde{M}$, which 
preserve $\widetilde{X}$, $\widetilde{\omega}$, and the lifted $C^\infty$ Lyapunov decomposition,
is a Lie group acting transitively on $\widetilde{M}$.}

{\it Proof.} By Proposition $2.7.$ of 
{\bf [BFL1]}, the $\nabla$-geodesics, tangent to $E^+$ or $E^-$, are complete, i.e. defined on $\mathbb{R}.$ 
Since $\nabla T=0$ and $\nabla R=0$, then by Lemma $A$ proved in the appendix, $\nabla$ is complete. 
So $\widetilde\nabla$ is also complete.

Recall that $E_a :=\{0\},$ if $a$ is not a Lyapunov exponent of $\phi_t.$ For $\forall\  a\in \mathbb{R}$, 
denote by $\widetilde P_a$ the projection of $T\widetilde M$ onto $\widetilde E_a.$ Since 
$\nabla E_a\subseteq E_a,$ then $\widetilde P_a$ is $\widetilde\nabla$-parallel. Thus 
$\{ \widetilde X,\  \widetilde\omega,\  \widetilde P_a\}_{a\in \mathbb{R}}$ is a family of 
$\widetilde\nabla$-parallel tensor fields. In addition, an application preserves 
$\{\widetilde P_a\}_{a\in \mathbb{R}}$, iff it preserves the lifted $C^\infty$ Lyapunov decomposition. 
So the lemma follows from 
the following classical result (see {\bf [K-No]}) :
 
{\it Let $N$ be a simply connected manifold, $\nabla_1$ be a complete 
connection on $N$ and $\mathcal{S}$ be a family of 
parallel tensor fields. If $\nabla_1 R^{\nabla_1}=0$ and 
$\nabla_1 T^{\nabla_1}=0$, then the group of $\nabla_1$-affine 
transformations which preserve $\mathcal{S}$ is 
a Lie group and acts transitively on $N$}. $\square$

$\ $

In the sense of the previous lemma, $\phi_t$ is called homogeneous. 
In particular, we deduce that $d\lambda$ vanishes nowhere. So on $M$, we have always 
two $C^\infty$ $\phi_t$-invariant line bundles $E^+_1$ and $E^-_1$, which are quite essential 
for the following discussions. \\\\
{\bf 4. The case of 2 positive Lyapunov exponents}\\
{\bf 4.1. Preparations}

Now we begin to eliminate the possible homogeneous models. In this section, we suppose that 
$\phi_t$ has two positive Lyapunov exponents. Then by Remark $3.2.1$, we have 
$$TM =\mathbb{R}X\oplus E^+_1\oplus E^+_2\oplus E^-_1\oplus E^-_2.$$
Up to a constant change of time scale, we suppose that the Lyapunov exponents of 
$E^+_1$ and $E^+_2$ are respectively $1$ and $a.$ 

In this case, the underlying geometric structure of our system is 
$$g_1 := (X, E^+_1, E^+_2, E^-_1, E^-_2, \omega).$$
Let $G'$ be the isometry group of $\widetilde g_1$ and $\Gamma$ be the 
fundamental group of $M$. By Lemma $3.2.4$, $G'$ acts 
transitively on $\widetilde M.$ The group $\Gamma$ is contained as a 
discrete subgroup in $G'$. Fix $x\in \widetilde{M}$ and 
denote by $H'$ the isotropy subgroup of $x$. Let $H'_e$ be the 
identity component of $H'$. Then we have the linear isotropy representation
$$ H'_e \stackrel{i}{\hookrightarrow} GL (T_x \widetilde{M})$$
$$h\longmapsto  D_x h.$$
Since each element of $H'$ preserves $\nabla$, then $i$ is injective. For $\forall\  
h\in H'_e$,
$$ D_xh(\widetilde X_x) =\widetilde X_x,\  D_xh(\widetilde E^\pm_x)\subseteq \widetilde E^\pm_x.$$
So in the following, we identify $i(h)$ with its restriction to $(\widetilde E^+\oplus \widetilde E^-)_x$.  

Take a basis $(l^+_2, l^+_1, l^-_2, l^-_1 )$ of $(\widetilde E^+\oplus \widetilde E^-)_x$, 
such that $l^\pm_{1, 2}\in (\widetilde{E}^\pm_{1, 2})_x$. Since each element $h$ 
of $H'_e$ preserves $\widetilde g_1$, then we have
$$ D_x h = \left (\begin{array}{cccc}
\lambda_1 & 0 & 0 & 0 \\
0 & \lambda_2 & 0 & 0 \\
0 & 0 & \frac{1}{\lambda_1} & 0\\
0 & 0 & 0 & \frac{1}{\lambda_2}
\end{array}
\right). $$
So $i( H'_e)$ is contained in a closed subgroup of $GL (T_x\widetilde{M})$,
which is isomorphic to $\mathbb{R}^2$. So we can identify $H'_e$ 
with $i(H'_e)$ and we deduce that $H'_e$ is isomorphic to $0$, 
$\mathbb{R}$ or $\mathbb{R}^2$. In any case, we have $\pi_1 (H'_e) = 0$.

Let $G'_e$ be the connected component of the identity of $G'$. Then it acts also 
transitively on $\widetilde{M}$. Using the long exact sequence of homotopy,
we get easily 
$$H'_e = H'\cap G'_e, \  \pi_1(G'_e) = 0.$$ 

Since $\widetilde{M}\cong G'/ H'$, then $\widetilde{M}$ is 
naturally equipped a real analytic structure. 
Since the geometric structure $\widetilde g_1$ is $G'$-invariant, then $\widetilde g_1$ is 
real analytic. Thus by {\bf [Am]} (see also {\bf [C-Q]}), the local Killing fields of $\widetilde g_1$ can 
be extended to global ones. Since $\nabla$ is in addition complete, then $H'$ is easily seen to have finitely 
many connected components. We deduce that $G'$ has also finitely many connected components. 
So up to finite covers, we can suppose that $\Gamma\subseteq G'_e.$

Denote by $\mathfrak{g}'$ and $\mathfrak{h}'$ the Lie algebras of $G'$ and $H'$. 
For $\forall\  u\in \mathfrak{g}'$, we have an induced $C^\infty$ Killing field on $\widetilde M$, 
$$Y^u : \widetilde M\to T\widetilde M$$
$$a\to \frac{d}{dt}\mid_{t=0}exp(tu)a.$$
Since $\nabla$ is complete, $\nabla R=0$ and $\nabla T=0$, 
then we have the following classical identification of vector spaces (see Theorem 2.8. of {\bf [K-No]} Ch.X)
$$j : \mathfrak{g}' \stackrel{\sim}{\longmapsto} T_x\widetilde{M}\oplus\mathfrak{h}'$$
$$ u \mapsto (Y^u(x), (\widetilde\nabla_{Y^u}-\mathcal{L}_{Y^u})\mid_x ),$$
where $\mathfrak{h}'$ has been identified with $Di(\mathfrak{h}')$ under $Di$. 

Pushing forward by $j$ the Lie algebra structure of $\mathfrak{g}'$ onto $T_x\widetilde{M}
\oplus\mathfrak{h}'$, we have for $\forall\  \  u,\  v \in T_x\widetilde{M}$ and 
$\forall\  \  A,\  B \in \mathfrak{h}'$, 
$$[ u, v ] = -T^{\widetilde\nabla}(u, v)- R^{\widetilde\nabla}(u, v), $$
$$[ A, u ] = A u, $$   
$$[A, B]= A\circ B - B\circ A.$$
 
Denote by $u$ the generating vector of the $1$-parameter subgroup $\{ \widetilde\phi_t\}_{t\in \mathbb{R}}$ 
of $G'$. Then $Y^u =\widetilde X$. Under the identification $j$, we have  
$$ u = \widetilde{X}_x + ( P^+_1- P^-_1+a P^+_2- a P^-_2) 
\in T_x\widetilde{M}\oplus \mathfrak{h}'.$$
If $L_0 := u -\widetilde{X}_x$, then $L_0\in \mathfrak{h}'.$ We deduce that $\mathfrak{h}'\cong \mathbb{R}$ or $\mathbb{R}^2.$

$\ $

{\bf Lemma 4.1.1.} {\it Under the notations above, $E^+_1\oplus E^-_1$ and 
$E^+_2\oplus E^-_2\oplus\mathbb{R}X$ are both integrable.}

{\it Proof.} Let $Y$, $Z$ be two $C^\infty$ sections of $E^+_1\oplus E^-_1$, then
$$ 0 = d\lambda (Y, Z) = -\lambda ([Y, Z]).$$
So $[Y, Z]$ is a section of $E^+\oplus E^-$.
$$i_{[Y, Z]}d\lambda = (\mathcal{L}_Y i_Z - i_Z\mathcal{L}_Y)d\lambda$$
$$= -i_Z(d i_Y + i_Y d) d\lambda = 0.$$
So $[Y, Z]$ is also a section of $E^+_1\oplus E^-_1$. 
Thus $E^+_1\oplus E^-_1$ is integrable.

Since $E^+_2$ and $E^-_2$ are both $\phi_t$-invariant, then $[X, E^\pm_2]\subseteq E^\pm_2.$ 
Define two tensor fields $K^\pm$ of type
$(1, 2)$ on $M$, such that 
$$K^\pm(Y, Z) = P^\pm_1[ P^+_2(Y), P^-_2(Z) ],\  \forall\  Y, Z\subseteq TM.$$
Then $K^\pm$ are both $\phi_t$-invariant. By Lemma $3.2.3$, 
$K^\pm( E^+_2, E^-_2 )
\subseteq \mathbb{R}X.$ 
So we have 
$$[E^+_2, E^-_2]\subseteq E^+_2\oplus E^-_2\oplus\mathbb{R}X.$$ 
Thus $E^+_2\oplus E^-_2\oplus\mathbb{R}X$ is integrable.  $\square$

$\ $

Up to finite covers, we suppose that $E^+$ and $E^-$ 
are both orientable. The connection $\nabla$ induces a connection $\nabla^+$ on $\wedge^2 E^+$. 
Denote by $\Omega^+$ its curvature form and by $\beta^+$  its connection form. Then we 
have 
$$\Omega^+( \cdot, \cdot )= {\rm Tr}( R( \cdot, \cdot)\mid_{E^+}),\  d\beta^+ = \Omega^+.$$

$\ $

{\bf Lemma 4.1.2.} {\it $d\lambda\wedge\Omega^+ = 0$, 
$\Omega^+\wedge\Omega^+ =0$, $\Omega^+\wedge\omega =0.$}

{\it Proof.} Since $\Omega^+$ is $\phi_t$-invariant and the flow is topologically
transitive, then there exists a constant $c$, such that
$$\lambda\wedge d\lambda\wedge\Omega^+ = c\cdot \lambda\wedge\omega\wedge\omega.$$
So 
$$c\int_M \lambda\wedge\omega\wedge\omega = \int_M \lambda\wedge d\lambda
\wedge\Omega^+$$
$$= -\int_M d(\lambda\wedge d\lambda\wedge\beta^+) =\int_{\partial M} 
\lambda\wedge d\lambda\wedge\beta^+ = 0.$$
So $c =0$. We deduce that $$d\lambda\wedge\Omega^+ = i_X(\lambda\wedge d\lambda
\wedge\Omega^+) =0.$$
 In the same way, we get $\Omega^+\wedge\Omega^+ =0$.\\
If $\lambda\wedge\Omega^+\wedge\omega = s \cdot\lambda\wedge\omega\wedge\omega$, then
$$ s\int_{M}\lambda\wedge\omega\wedge\omega = \int_M \beta^+\wedge 
d\lambda\wedge\omega.$$
If $\lambda\wedge d\lambda\wedge\omega =\delta \cdot\lambda\wedge\omega\wedge\omega,$
then
$$\beta^+\wedge d\lambda\wedge \omega = \delta \cdot\beta^+\wedge\omega\wedge\omega$$
$$=\delta\cdot \beta^+(X)\lambda\wedge\omega\wedge\omega. $$
By the same argument as in Lemma $2.3.3.$ of {\bf [BFL2]}, we get 
$$\int_{M}\beta^+(X)\lambda\wedge\omega\wedge\omega =0.$$
So $s =0$, i.e. $\Omega^+\wedge\omega =0.$  $\square$

$\ $

{\bf Lemma 4.1.3.}  {\it Under the notations above, we have $\Omega^+ =0.$}

{\it Proof.} In the direction of $X$, the situation is always clear. So in the following, we consider 
only the restrictions onto $E^+\oplus E^-$ of the forms and endomorphisms. 

Since $\omega\mid_{E^+\oplus E^-}$ is non-degenerate, then we can 
find a section $\psi$ of $End(E^+\oplus E^-),$ such that
$$\Omega^+(\cdot, \cdot) = \omega (\psi(\cdot),\  \cdot).$$
For $\forall\  y\in M$, take $l^\pm_{1, 2}\in (E^\pm_{1, 2})_y$ such 
that  $(l^+_2, l^+_1, l^-_2, l^-_1)$ forms a dual basis of $\omega_y$, i.e. 
$$\omega(l^+_2, l^-_2) =\omega(l^+_1, l^-_1) =1,\  \omega(l^+_2, l^-_1) =\omega(l^+_1, l^-_2) =0.$$

If $\psi_y(l^+_1) =0,$ then in this basis, we get
$${\psi_y}=\left (
\begin{array}{cccc}
A & 0 & 0 & 0\\
B & 0 & 0 & 0\\
0 & 0 & A & B\\
0 & 0 & 0 & 0
\end{array}\right )$$
Since $\Omega^+\wedge\omega =0,$ then Tr$\psi = 2A=0.$ 
By Lemma $2.1.3$, 
$$ 0 = \Omega^+_y(l^+_2, l^-_1) = \omega (\psi l^+_2, l^-_1)$$
$$=B \cdot \omega(l^+_1, l^-_1).$$
So $B =0$. Thus $\psi_y =0.$

Now suppose that $\psi_y(l^+_1)\not =0$. Since $\Omega^+\wedge\Omega^+ =0$,
then det$(\psi_y) =0$. So $\exists\  \  y^+_1= \alpha l^+_2 +\delta l^+_1,$ $\alpha\not =0$, 
such that $\psi_y( y^+_1) =0$. Then in a dual basis with respect to $\omega_y$, 
$(y^+_1, l^+_1, y^-_1, z^-),$ we have 
$${\psi_y}=\left(
\begin{array}{cccc}
0 & A & 0 & 0\\
0 & B & 0 & 0 \\
0 & 0 & 0 & 0 \\
0 & 0 & A & B
\end{array}\right)$$
As above, we have Tr$(\psi_y) =2B =0$. By Lemma $2.1.3$, 
$$0 =\Omega^+_y(l^+_1, l^-_2) =\omega(A y^+_1, l^-_2)$$
$$=A\cdot\alpha\cdot\omega(l^+_2, l^-_2) =A\cdot\alpha.$$ 
So $A =0.$ We deduce that $\psi  \equiv0$, i.e. $\Omega^+\equiv 0.$  $\square$

$\ $

Define the following map
$$\mathfrak{g}'\stackrel{\chi}{\longmapsto}\mathbb{R}$$
$$ u+A \mapsto Tr(A\mid_{\widetilde{E}^+_x}).$$
Since $\Omega^+ \equiv 0$, then $\chi$ 
is a character of $\mathfrak{g}'$. So the kernel of $\chi$ is an ideal of $\mathfrak{g}'$, 
denoted by $\mathfrak{g}$, 

We have seen that $\mathfrak{h}'$ is isomorphic to $\mathbb{R}$ or 
$\mathbb{R}^2$. In the following, these two cases are considered 
seperately.\\\\
{\bf 4.2. dim$\mathfrak{h}'$ =1.}

In this subsection, we suppose that dim$\mathfrak{h}' =1.$ To prove 
the non-existence of such a flow, we shall at first calculate explicitly $\mathfrak{g}'$ using the lemmas 
estalished in the previous subsection. Then we shall get a contradiction via the non-existence 
of cocompact lattice in $\mathbb{R}^2\rtimes \widetilde{SL(2, \mathbb{R})}.$
    
Since $L_0\in \mathfrak{h}'$ (see Subsection $4.1$), 
then $\mathfrak{h}'=\mathbb{R}L_0.$ To simplify the notations, we identify $T_x\widetilde M$ with $T_xM$. 
Thus we have 
$$\mathfrak{g}' =T_xM\oplus \mathfrak{h}'.$$
Denote by $\mathfrak{g}$ the kernel of $\chi$ (see Subsection $4.1$). Then $\mathfrak{g}$ is an ideal of 
$\mathfrak{g}'$. Since $\chi(L_0) =1+a >0$, then we have $\mathfrak{g} =T_xM.$ Recall that the Lyapunov 
exponents of $E^+_1$ and $E^+_2$ are $1$ and $a$. Now we can find explicitly $\mathfrak{g}$ as follows.

Since $\mathfrak{g}$ ($= T_xM$) is an ideal of $\mathfrak{g}'$, then for $\forall\  u,\  v\in T_xM$, 
$$[u, v] = -T(u, v) -R(u, v)\in T_xM.$$
Thus $R(u, v) =0$ and $[u, v]= -T(u, v)$.

Take a basis of $T_xM$, $(X_x, l^+_2, l^+_1, l^-_2, l^-_1)$, such that $l^\pm_{1, 2}\in (E^\pm_{1, 2})_x$ 
and $d\lambda (l^+_2, l^-_2) = -1$. 
Extend $l^\pm_{1, 2}$ to local sections of $E^\pm_{1, 2}$, 
denoted by ${\bar l}^\pm_{1, 2}$. By the definition of $\nabla$, we get 
$$[X_x, l^\pm_1 ] = - T(X_x, l^\pm_1) =\mp l^\pm_1.$$
Similarly, 
$$[X_x, l^\pm_2 ]= \mp a l^\pm_2.$$
Since $E^+_1\oplus E^-_1$ is integrable by Lemma $4.1.1$, then we get
$$[l^+_1, l^-_1] = -T(l^+_1, l^-_1)$$
$$= -( P^-_1[{\bar l}^+_1, {\bar l}^-_1] + P^+_1[{\bar l}^+_1, {\bar l}^-_1]- [{\bar l}^+_1, {\bar l}^-_1]) =0.$$
Similarly, we get
$$[l^+_2, l^-_2] = X_x.$$

$\ $

{\bf Lemma 4.2.1.} {\it Under the notations above, we have $1< a$.}

{\it Proof.} Suppose that $1>a.$ Then by Lemma $3.2.3$, we have 
$$T(E^+_1, E^-_2)\subseteq E_{1-a}.$$
If $1-a\not =a$, then $[l^+_1, l^-_2]= -T(l^+_1, l^-_2) =0.$ If $1-a =a$, then 
$\exists \  b\in \mathbb{R},$ such that 
$$[l^+_1, l^-_2] =b\cdot l^+_2.$$
So in any case, $\exists \  c\in \mathbb{R}$, such that $[l^+_1, l^-_2] =c\cdot l^+_2.$

Since $T(E^+_1, E^+_2)\subseteq E_{1+a} =\{0\},$ then $[l^+_1, l^+_2] =0.$ By the Jacobi identity of 
$l^+_1$, $l^+_2$ and $l^-_2$, we get 
$$0=[l^+_2, [l^+_1, l^-_2]] + [l^+_1, [l^-_2, l^+_2]] +[l^-_2, [l^+_2, l^+_1]]$$
$$=[l^+_2, c\cdot l^+_2] + [l^+_1, - X_x] = - l^+_1,$$
which is absurd.  $\square$
 
$\ $

Since $a >1$, then we can suppose that $[l^+_1, l^-_2] =c\cdot l^-_1$ and $[l^-_1, l^+_2] =d\cdot l^+_1.$ 
Again by the Jacobi identity of $l^+_1$, $l^+_2$ and $l^-_2$, we have 
$$0 =[l^+_2, [l^+_1, l^-_2]] +[l^+_1, [l^-_2, l^+_2]] = -(1 +c\cdot d) l^+_1.$$
So $c\cdot d = -1$. Now replacing $l^-_2$ by $\frac{1}{c}l^-_2$ and $l^+_2$ by $c\cdot l^+_2$, 
we get the following bracket relations of $\mathfrak{g}$, 
$$[X_x, l^\pm_1] =\mp  l^\pm_1,\   [X_x, l^\pm_2] =\mp a l^\pm_2,$$
$$[l^+_1, l^-_1] =0,\   [l^+_1, l^-_2] = l^-_1,$$
$$[l^-_1, l^+_2] = - l^+_1,\   [l^+_2, l^-_2] = X_x.$$
The brackets, not appeared above, vanish by Lemma $3.2.3.$

Since $[l^+_1, l^-_2] =l^-_1$, then 
$E_{1-a}\not =\{0\}.$ We deduce that $a=2.$ Thus by the bracket 
relations above, we get clearly
$$\mathfrak{g}\cong \mathbb{R}^2\rtimes\mathfrak{sl}(2, \mathbb{R}),$$
where the semi-direct product is given by matrix multiplication.

It is easily seen that the center of 
$\mathfrak{g}'$ is $\mathbb{R} (X_x+ L_0)$. Thus we have the following 
direct product decomposition
$$\mathfrak{g}' \cong \mathfrak{g}\oplus \mathbb{R}(X_x + L_0).$$ 
Let $G$ be the connected subgroup of $G'_e$ integrating $\mathfrak{g}$. 
Since $G'_e$ is simply connected (see Subsection $4.1$), 
then $G$ is also simply connected and $G'_e = G\times \mathbb{R}$, where $\mathbb{R}$ integrates 
$\mathbb{R}(X_x +L_0)$ in $G'_e$. Thus we get  
$$G \cong \mathbb{R}^2\rtimes\widetilde{SL(2, \mathbb{R})}.$$
It is easily seen that $G$ acts 
transitively on $\widetilde M$. Then by the long exact sequence of homotopy, $G\cap H'_e$ is 
seen to be connected. 
So $G\cap H'_e = \{ e\},$ i.e. $G$ acts freely on $\widetilde M$. Thus $G$ is identified to $\widetilde{M}.$ 

Up to finite covers, we have $\Gamma\subseteq G'_e$ (see Subsection $4.1$). Let $\Gamma_1$ 
be the projection of $\Gamma$ into $G$, with respect to the direct product $G'_e = G\times \mathbb{R}$. 
Since $\Omega^+ =0$, 
then by the general arguments of Section $5.$ of {\bf [BFL1]}, 
$\Gamma_1$ is seen to be a cocompact lattice of $G$. Now we eliminate 
this case by proving

$\ $

{\bf Lemma 4.2.2.} {\it $\mathbb{R}^2\rtimes\widetilde{SL(2,\mathbb{R})}$ 
has no cocompact lattice.}

{\it Proof.} Suppose that there exists a cocompact lattice, denoted by $\Delta.$ 
Define $\Delta_1:=\Delta\cap \mathbb{R}^2$ and denote by $\Delta_2$ the projection 
of $\Delta$ to $\widetilde{SL(2,\mathbb{R})}$. 
Then by Corollary $8.28$ of {\bf [Ra]}, 
$\Delta_1$ is a cocompact lattice of $\mathbb{R}^2$ 
and $\Delta_2$ is a lattice of $\widetilde{SL(2,\mathbb{R})}$. 

Denote by $\pi$ the natural projection of $\widetilde{SL(2,\mathbb{R})}$ onto $SL(2,\mathbb{R})$. Then 
$\pi(\Delta_2)$ preserves the lattice $\Delta_1$ for the linear action. So $\pi(\Delta_2)$ is conjugate to a 
subgroup of $SL(2, \mathbb{Z})$.

Since $\Delta$ is cocompact, then $\pi(\Delta_2)$ is also cocompact. We deduce that 
$SL(2, \mathbb{Z})$ is cocompact in $SL(2, \mathbb{R})$, which is absurd.  $\square$\\\\
{\bf 4.3. dim$\mathfrak{h}'$ =2.}

In this subsection, we suppose that dim$\mathfrak{h}' =2.$ To prove the non-existence of such a flow, 
we shall at first find $\mathfrak{g}'$. Then we shall study the action of the fundamental 
group of $M$ on the space of lifted weak unstable leaves to deduce a dynamical contradiction. 

Define $S : = P^+_2 -P^+_1 - P^-_2 +P^-_1$. Then 
$\mathfrak{h}'$ is generated by $S$ and $L_0$(see Subsection $4.1.$). Since 
$\chi(S )=0$, then we have $\mathfrak{g} = \mathbb{R}S 
\oplus T_x M$.

As in the previous subsection, we suppose that the Lyapunov 
exponents of $E^+_1$ and $E^+_2$ are $1$ and $a$. Take a basis of $T_xM$, $(X_x, l^+_2, l^+_1, l^-_2, l^-_1),$ 
such that $ l^\pm_{1, 2}\in 
E^\pm_{1, 2}$ and $d\lambda (l^+_2, l^-_2) = -1.$ Suppose at first that $ a>1.$ Then by the same argument as 
in Lemma $4.2.1$, we can find $c$ and $d$, such that 
$$ [l^+_1, l^-_2] = c\cdot l^-_1,\  [l^-_1, l^+_2] =d\cdot l^+_1.$$
By the Jacobi identity of $S$, $l^+_1$ and $l^-_2$, we get 
$$0=[S, [l^+_1, l^-_2]] +[l^+_1, [l^-_2, S]] +[l^-_2, [S, l^+_1]]$$
$$=c\cdot l^-_1 +[l^+_1, l^-_2] +[l^-_2, -l^+_1] = 3c\cdot l^-_1.$$
Thus $c=0.$ Similarly we get $d=0$. If $a<1$, then we can find $c'$ and $d'$, such that 
$$ [l^+_1, l^-_2] = c'\cdot l^+_2,\  [l^-_1, l^+_2] =d'\cdot l^-_2.$$
Thus by the Jacobi identities, we get as above $c'=d'=0$. We deduce that 
$$  [l^+_1, l^-_2] =0,\  [l^-_1, l^+_2] =0.$$

Now by similar arguments as in the previous subsection, we get the 
following bracket relations, 
$$ [S, l^\pm_1] = \mp l^\pm_1,\   [S, l^\pm_2] =\pm l^\pm_2,$$
$$[L_0, l^\pm_{1}] =\pm  l^\pm_1,\   [L_0, l^\pm_2] =\pm a l^\pm_2,$$
$$[X_x, l^\pm_1] =\mp l^\pm_1,\   [X_x, l^\pm_2] =\mp a l^\pm_2,\  [l^+_2, l^-_2] = X_x- S.$$
The brackets, not appeared above, vanish. Define three elements,
$$ \alpha :=\frac{L_0 + S}{a + 1},\  \beta :=\frac{L_0- a S}{a+1}, \  \delta :=\frac{X_x - S}{a+1}.$$ 
Then $\mathfrak{g}'$ is decomposed 
as a direct product of three ideals,
$$ \mathfrak{g}'\cong (\mathbb{R}l^+_1\oplus \mathbb{R} l^-_1\oplus
\mathbb{R}\beta)\oplus (\mathbb{R}l^+_2\oplus\mathbb{R}l^-_2
\oplus\mathbb{R}\delta)\oplus \mathbb{R}(\delta +\alpha). $$
Then by the bracket relations above, we get 
$$ \mathfrak{g}'\cong (\mathbb{R}^2\rtimes \mathbb{R})\times \mathfrak{sl}(2, \mathbb{R})\times\mathbb{R},$$
where the semi-direct product, $\mathbb{R}^2\rtimes\mathbb{R}$, is 
given by the linear action on $\mathbb{R}^2$ 
of the order two diagonal matrices of trace zero. 
Since $G'_e$ is simply-connected, then we have 
$$ G'_e\cong (\mathbb{R}^2\rtimes \mathbb{R})\times
\widetilde{SL(2, \mathbb{R})}\times\mathbb{R}.$$

Now we begin to study the action of $\Gamma$ on the space of lifted weak unstable leaves. Let us recall at first 
some notations.

Let $\psi_t$ be a $C^\infty$ Anosov flow on a closed manifold $N$. Denote by 
$\widetilde\psi_t$ its lifted flow on the universal covering space $\widetilde N.$ 
Denote by $\widetilde{\mathcal{F}}^{+, 0}$ the lifted foilation of $\mathcal{F}^{+, 0}$ 
and by ${\widetilde N}/{\widetilde{\mathcal{F}}^{+, 0}}$ the space of lifted weak unstable leaves with the quotient 
topology. Thus the fundamental group $\pi_1(N)$ acts naturally on ${\widetilde N}/{\widetilde{\mathcal{F}}^{+, 0}}$. 
The following lemma has appeared in some special 
contexts (see for example {\bf [BFL2]} and {\bf [Ba]}). For the sake of completeness, we prove it in detail.

$\ $

{\bf Lemma 4.3.1.} {\it Under the notations above, if 
$\gamma\in\pi_1(N)$ and $\gamma\not=e$, then each $\gamma$-fixed point 
of ${\widetilde N}/{\widetilde{\mathcal{F}}^{+, 0}}$ is either contractive or repulsive.}

{\it Proof.} Suppose that $\widetilde{W}^{+, 0}_x$ is fixed by $\gamma.$ Then $\exists\  t\in\mathbb{R}$, such that 
$$\gamma \widetilde{W}^{+}_x = \widetilde\phi_t  \widetilde{W}^{+}_x.$$
If $t =0$, then we can take a curve $l$ in $\widetilde{W}^{+}_x$, such that $l(0)=x$ and $l(1) =\gamma x.$ If 
$s<<0$, then $\phi_s(\pi(l))$ will be tiny, where $\pi$ denotes the projection of $\widetilde N$ onto $N$. Thus 
$\phi_s(\pi(l))$ is homotopically trivial. We deduce that $\pi(l)$ is also homotopically trivial, i.e. $\gamma=e$, which 
is a contradiction. So $t\not =0.$

By replacing $\gamma$ by $\gamma^{-1}$ if necessary, we suppose that $t<0.$ We can see as follows that 
$\widetilde{W}^{+, 0}_x$ is $\gamma$-contractive.

Fix a $C^\infty$ Riemannian metric $g$ on $N$. Denote by $\widetilde g$ the lifted metric on $\widetilde N$. By 
{\bf [An]}, the induced metrics on the leaves of $\mathcal{F}^{+, 0}$ are all complete. Thus with its 
induced metric, $\widetilde{W}^{+}_x$ is a complete metric space. Since $\gamma$ acts isometrically, then 
$\gamma^{-n}\circ\widetilde\phi_{nt}$ is a contraction of $\widetilde{W}^{+}_x$, if $n>>1.$ Thus it admits a 
unique fixed point in $\widetilde{W}^{+}_x$, denoted again by $x$. So we get 
$$ \gamma x = \widetilde\phi_t  x,$$
i.e. the orbit of $x$ is fixed by $\gamma.$ 

Denote by $\bar U$ the saturated set of $\widetilde{W}^{-}_x$ with respect to $\widetilde{\mathcal{F}}^{+, 0}$. 
Then by 
the local product structure of $\widetilde\phi_t$, $\bar U$ is open. Thus the 
projection of $\widetilde{W}^{-}_x$ into $\widetilde N/\mathcal{F}^{+, 0}$ is an open 
neighborhood of $\widetilde{W}^{+, 0}_x$, denoted by $U$. For $\forall\   y\in \widetilde{W}^-_x,$ we have 
$$\gamma^n \widetilde{W}^{+, 0}_y = \widetilde{W}^{+, 0}_{(\widetilde\phi_{-nt}\circ\gamma^n)(y)}.$$
Since $ \gamma x = \widetilde\phi_t  x,$, then $(\widetilde\phi_{-nt}\circ\gamma^n)(x)=x.$ 
So $(\widetilde\phi_{-nt}\circ\gamma^n)(y)\xrightarrow[n\rightarrow +\infty]{} x.$ 
We deduce that $\gamma^n \widetilde{W}^{+, 0}_y \xrightarrow[n\rightarrow +\infty]{} \widetilde{W}^{+,0}_x.$ 
So $\gamma$ contracts on $U$.  $\square$

$\ $

Now return to our {\it geometric} Anosov flow $\phi_t$. Denote by $P'_e$ the stabilizer of $\widetilde{W}^{+}_x$ in 
$G'_e$. Then $H'_e\subseteq P'_e$ and $P'_e$ is easily seen to be connected. So 
$G'_e/P'_e$ is identified to $\widetilde M/\widetilde{\mathcal{F}}^{+, 0}$. Define 
$$\mathfrak{p}^+
:= \mathbb{R}X_x\oplus \mathfrak{h}'\oplus\mathbb{R}l^+_1
\oplus\mathbb{R}l^+_2.$$
Then $\mathfrak{p}^+$ is the Lie algebra of $P'_e$ and $P'_e$ is seen to be closed 
in $G'_e$. Since $G'_e$ is simply connected (see Subsection $4.1$), then by the long exact sequence 
of homotopy, we get $\pi_1(G'_e/P'_e) =0.$ 

Define $G^1_e := (\mathbb{R}^2\rtimes\mathbb{R})\times 
SL(2, \mathbb{R})\times \mathbb{R}$ and 
denote by $P^1_e$ the connected Lie subgroup of $G^1_e$ with Lie algebra 
$\mathfrak{p}^+$. Then $G^1_e/P^1_e$ is naturally identified to $\mathbb{R}^1\times \mathbb{S}^1$. Denote by 
$\pi$ the projection of $G'_e$ onto $G^1_e$ and by $P''_e$ the group $\pi^{-1}(P^1_e)$. Then we get 
$$  G'_e/P''_e\cong G^1_e/P^1_e\cong  \mathbb{R}^1\times \mathbb{S}^1.$$
We deduce that $G'_e/P'_e\cong   \mathbb{R}^1\times \mathbb{R}^1.$

Since $\phi_t$ preserves a volume form, then the periodic points of $\phi_t$ is dense in $M$. 
Take $g H'_e\in G'_e/H'_e$ ($\cong \widetilde M$), such that its projection in $M$ is of period $T$. 
If $\widetilde\phi_T(gH'_e) =gH'_e$, then each orbit of $\widetilde\phi_t$ is periodic by the 
homogeneity of $\widetilde \phi_t$. We deduce that 
each $\phi_t$-orbit is periodic, which contradicts the topological transitivity of $\phi_t$. So 
$\widetilde\phi_T(gH'_e) \not =gH'_e.$

Now take $\gamma\in \Gamma$ $(\subseteq G'_e)$, such that $\gamma(gH'_e)=\widetilde\phi_T(gH'_e)$. 
Then we have $\gamma\not =e$ and $\exists\  h\in H'_e$, such that 
$$ \gamma =g(h\cdot exp(T(X_x+L_0))g^{-1}.$$
Since $\gamma$ fixes the orbit of $gH'_e$, then it fixes $gP'_e$ and $gP''_e$. So by Lemma $4.3.1$, 
the $\gamma$-action on $G'_e/P''_e$ admits at least an isolated fixed point. Then 
by some direct calculations, the corresponding $\gamma$-action on 
$\mathbb{R}^1\times \mathbb{S}^1$ ($\cong G'_e/P''_e$) must be as following 
$$ \mathbb{R}^1\times \mathbb{S}^1\stackrel{\gamma}{\to} \mathbb{R}^1\times \mathbb{S}^1 \eqno(\ast)$$
$$(y,\  [u])\to ( exp^{-c}y+d, \  [Au]),$$
where $c\not =0$ and $A$ is a matrix with two different positive eigenvalues. Here $\mathbb{S}^1$ is viewed 
as the set of directions, i.e.
$$  \mathbb{S}^1\cong \{ [u]\mid u\in \mathbb{C}^\ast, u\sim v\Leftrightarrow u = t v, t>0\}.$$
Then $GL(2, \mathbb{R})$ acts on $\mathbb{S}^1$ by matrix multiplication. 

Up to an isomorphism of covering spaces, the projection of $G'_e/P'_e$ onto $G'_e/P''_e$ 
is as following  
$$\mathbb{R}^1\times\mathbb{R}^1\longmapsto\mathbb{R}^1
\times\mathbb{S}^1 \eqno(\ast\ast)$$
$$( x, \  \theta )\mapsto ( x, \  [e^{i\theta}]).$$
Since the $\gamma$-action on $G'_e/P'_e$ is just a lift of the $\gamma$-action on $G'_e/P''_e$, then by 
$(\ast)$ and $(\ast\ast)$, we clearly see that on $G'_e/P'_e$, $\gamma$ admits either a saddle or no fixed point. 
We deduce that $\gamma$ admits a saddle on $G'_e/P'_e$, which contradicts Lemma $4.3.1.$\\\\ 
{\bf 5. The case of 1 positive exponent and $d\lambda\wedge\omega\not\equiv 0$.}\\
{\bf 5.1. Preparations.}

In this section, we suppose that $\phi_t$ has only one positive Lyapunov exponent and 
$d\lambda\wedge\omega\not\equiv 0$. Up to 
a constant change of time scale, we suppose that 
this positive exponent is $1$. By Lemma $3.2.4$, $d\lambda\wedge\omega$ vanishes 
nowhere. So $\omega\mid_{E^+_1\oplus E^-_1}$ is non-degenerate. As in Lemma $3.2.2$, we define 
$$ E^\pm_2 :=[ Ker(v\mapsto\omega(E^\mp_1, v))]\cap E^\pm.$$
Then $E^+_2$ and $E^-_2$ are both $\phi_t$-invariant $C^\infty$ line subbundles of $TM.$ 

In this case, the underlying geometric structure is 
$$g_2 := (X, E^+, E^-, \omega).$$
Denote by $G'$ the isometry group of $\widetilde g_2$. Then by Lemma $3.2.4$, $G'$ acts 
transitively on $\widetilde M.$ Fix $x\in \widetilde M$ and denote by $H'$ the isotropy subgroup of 
$x$. Because of the existence of $E^\pm_2$, some arguments of Subsection $4.1.$ pass through 
without change. In particular, we get that $H'_e$ is isomorphic to $0$, $\mathbb{R}$ or $\mathbb{R}^2$ and 
$G'_e$ is simply connected.

Denote by $\mathfrak{g}'$ and $\mathfrak{h}'$ the Lie algebras of $G'$ and $H'$. By using the 
connection corresponding to the case of one positive Lyapunov exponent, we have a similar 
identification of $\mathfrak{g}'$ and $T_x\widetilde M\oplus\mathfrak{h}'$ as in Subsection $4.1$. To simplify 
the notations, we identify $T_x\widetilde M$ with $T_xM.$ If $L_0 : =P^+ -P^-$, 
then we get $L_0\in\mathfrak{h}'$. So $\mathfrak{h}'\cong \mathbb{R}$ or $\mathbb{R}^2$.

Lemmas $4.1.1.$ and $4.1.2.$ are also valid here. But the proof of Lemma $4.1.3.$ 
does not pass through in the current case.\\\\
{\bf 5.2. dim$\mathfrak{h}'$ =2.}

In this subsection, we suppose that $\mathfrak{h}'$ is of dimension two. So if we 
define $S :=P^+_2- P^+_1-P^-_2 +P^-_1$, then $\mathfrak{h}'$ is generated by $S$ and $L_0.$

$\ $

{\bf Lemma 5.2.1.} {\it Under the notations above, we have $\Omega^+\equiv 0.$}

{\it Proof.} As Before, we consider only the restrictions onto $E^+\oplus E^-$ of the forms and endomorphisms. Take 
a dual basis with respect to $\omega_x\mid_{{E}^+_x
\oplus{E}^+_x}$, $(l^+_2, l^+_1, l^-_2, l^-_1 )$, such that 
$l^\pm_{1, 2}\in{E}^\pm_{1, 2}$ and $d\lambda(l^+_2, l^-_2) =-1.$ Extend locally these 
vectors to the sections of $E^\pm_{1, 2}$, denoted by $\bar l^\pm_{1, 2}.$ Then we have 
$$ T(l^+_1, l^-_2) = P^-[\bar l^+_1, \bar l^-_2]-P^+[ \bar l^-_2, \bar l^+_1] -[\bar l^+_1, \bar l^-_2]$$
$$=d\lambda(l^+_1, l^-_2 )\cdot X_x=0.$$
Similarly, we have $T(l^+_1, l^-_1)=0$ and $T(l^+_2, l^-_2) = -X_x.$ Thus we get the constants, such that 
$$[l^+_1, l^-_1]=aS +b L_0, \  [l^+_1, l^-_2]=a'S +b' L_0,$$
$$[ l^+_2,  l^-_2]= X_x +a''S +b'' L_0.$$
As before, we have 
$$[l^\pm_1, l^\pm_2] =0,\  [X_x, l^\pm_{1, 2}]=\mp l^\pm_{1, 2}.$$
By the Jacobi identity of $l^+_1$, $l^+_2$ and $l^-_2$, we get
$$ 0=[l^+_1, [l^+_2, l^-_2]]+[l^+_2, [l^-_2, l^+_1]] +[l^-_2, [l^+_1, l^+_2]]$$
$$=[l^+_1, X_x +a''S+ b''L_0]+[l^+_2, -a'S-b'L_0]$$
$$= (1+a''-b'')l^+_1 +(a'+b')l^+_2.$$
So $a'+b'=0.$ By the Jacobi identity of $l^-_2$, $l^+_1$ and $l^-_1$, 
we get $a'-b'=0.$ So $[l^+_1, l^-_2] =0.$ We deduce that $\Omega^+(l^+_1, l^-_2) =0.$ 

Define $\psi$ as in Lemma $4.1.3$. View $\psi_x$ as a matrix in 
the basis above, then we get $(\psi_x)_{1,2} = \Omega^+(l^+_1 , l^-_2)=0.$ Since 
$\Omega^+\wedge\Omega^+ = \Omega^+\wedge\omega=0$, then det$\psi$ =Tr$\psi$ =$0.$ 
So we get   
$$\psi_x =\left (\begin{array}{cccc}
0 & 0 & 0 & 0\\
B & 0 & 0 & 0\\
0 & 0 & 0 & B\\
0 & 0 & 0 & 0
\end{array}\right )$$
For $\forall\  h\in H'_e$, $h$ preserves $\Omega^+$. If $\psi_x\not =0,$ 
then the matrix of $D_xh$ must have the following form
$$D_xh=\left (\begin{array}{cccc}
c & 0 & 0 & 0 \\
d & c & 0 & 0\\
0 & 0 & \frac{1}{c} & -\frac{d}{c^2} \\
0 & 0 & 0 & \frac{1}{c} \end{array}\right ) $$
But $h$ preserves also the subbundles, $E^\pm_{1, 2}$. So $d=0$. 
We deduce that dim$\mathfrak{h}' = 1$, which is a contradiction.
So we get $\psi_x =0,$ i.e. $\Omega^+_x=0.$ Then 
by homogeneity, $\Omega^+ \equiv 0.$   $\square$

$\ $

With the help of the previous lemma, we can define as in Subsection $4.1.$ a 
character $\chi$ of $\mathfrak{g}'$. Then by 
similar calculations as in Subsection $4.2$, $\mathfrak{g}'$ is seen to be the same as that 
of Subsection $4.3$, except that $a=1$ here. But in Subsection $4.3$, we have 
found three elements, $\alpha$, $\beta$ and $\delta$, which have eliminated the effect of $a$ on 
the structure of $\mathfrak{g}'$. So we get here the same $G'_e$ and $H'_e$ as in Subsection $4.3$. 
Thus the same arguments prove the non-existence of this case.\\\\
{\bf 5.3. dim$\mathfrak{h}'$ =1.}

In this subsection, we suppose that dim$\mathfrak{h}' =1.$ So $\mathfrak{g}' =\mathbb{R}L_0 \oplus T_xM.$ 
Take a basis $(X_x, l^+_2, l^+_1, l^-_2, l^-_1)$ of $T_xM$, such that 
$l^\pm_{1, 2}\in (E^\pm_{1, 2})_x$ and $d\lambda(l^+_2, l^-_2) =1$. Since 
$E^+_2$ and $E^-_2$ are both of dimension $1$, then there exists a well-defined 
smooth function $f$, such that 
$$d\lambda\mid_{E^+_2\oplus E^-_2} =f\cdot\omega\mid_{E^+_2\oplus E^-_2}.$$ 
Since $d\lambda$, $\omega$, and $E^\pm_2$ are all $\phi_t$-invariant, 
then $f$ is also $\phi_t$-invariant. We deduce that $f$ is constant. So if we multiply $\omega$
by a constant, $(l^+_2, l^+_1, l^-_2, l^-_1)$ can be supposed to
be dual with respect to $\omega_x\mid_{E^+_x\oplus E^-_x}$. 
Using the Jacobi identities, we get directly (see Subsection $4.2.$)
$$[ L_0, l^\pm_{1, 2}] = \pm l^\pm_{1, 2},$$
$$[X_x, l^\pm_{1, 2}] =\mp l^\pm_{1, 2},\   [l^+_2, l^-_2] = -X_x - L_0.$$
The brackets, not appeared above, vanish. 

Define $\alpha :=X_x + L_0$ and $\mathfrak{g}\cong \mathbb{R}\alpha\oplus\mathbb{R} l^+_2\oplus
\mathbb{R} l^+_1\oplus\mathbb{R} l^-_2\oplus\mathbb{R} l^-_1$. Then we get
$$\mathfrak{g}'\cong \mathfrak{g}\rtimes \mathbb{R} L_0.$$ 
Denote by $G_e$ the connected Lie subgroup of $G'_e$ with Lie algebra $\mathfrak{g}$. Since 
$G'_e$ is simply connected, then so is $G_e$. Thus by the bracket relations above, we get
$$G_e\cong \mathbb{R}^2\times Heis,$$
where $Heis$ represents the three dimensional Heisenberg group. In addition, 
we have $G'_e\cong G_e\rtimes H'_e$. So $G_e$ is naturally 
identified to $\widetilde M$ as follows
$$ \psi : G_e\stackrel{\sim}{\longmapsto} \widetilde M$$
$$g \mapsto gx.$$

Define $\omega_1:= \psi^\ast\omega$. Then $\omega_1$ is a left-invariant $2$-form on $G_e$. 
View $l^\pm_{1, 2}$ as left-invariant vector fields on $G_e$. Then 
$\mathbb{R}l^+_2\oplus \mathbb{R} l^+_1$ and $\mathbb{R} l^-_2\oplus\mathbb{R}l^-_1$ 
are identified to $\widetilde E^+$ and $\widetilde E^-$. The corresponding flow on $G_e$ is given by 
the left-invariant field $\alpha.$ 
So the corresponding geometric structure on $G_e$ is given by
$$g_3 :=( \alpha, \mathbb{R}l^+_1\oplus\mathbb{R}l^+_2, \mathbb{R}l^-_1\oplus
\mathbb{R}l^-_2, \omega_1).$$
In addition, by the identification of 
$\mathfrak{g}'$ with $\mathfrak{h}'\oplus T_xM$, $(l^+_2, l^+_1, l^-_2, l^-_1)$ 
is dual with respect to $\omega_1$

For $\forall\  c, d\in\mathbb{R}$, define an endomorphism of $\mathfrak{g}$, $\rho^c_d$, such that 
$$ \rho^c_d( l^\pm_1) = exp^{\pm c} l^\pm_1,$$
$$\rho^c_d(l^\pm_2) =  exp^{\pm d} l^\pm_2, \  \rho^c_d (\alpha) = \alpha .$$
Then $\{\rho^c_d\}_{c, d\in \mathbb{R}}$ gives a two-parameter family of Lie 
algebra automorphisms of $\mathfrak{g}$. The corresponding isomorphisms of $G_e$ 
form a Lie group isomorphic to $\mathbb{R}^2$. Then 
we observe easily that its action on $G_e$ preserves $g_3$ and fixes $e$. We deduce that 
dim$(H'_e) \geq 2$, which is contradictory to the assumption that dim$\mathfrak{h}' =1$.\\\\ 
{\bf 6. The case of 1 positive exponent and $d\lambda\wedge\omega\equiv 0.$}\\
{\bf 6.1. Preparations.}

In this section, we suppose that $\phi_t$ has only one positive Lyapunov exponent and 
$d\lambda\wedge\omega\equiv 0.$ As before, we suppose that this positve exponent is $1.$ 
Since $d\lambda\wedge\omega\equiv 0,$ then $\omega\mid_{E^+_1\oplus E^-_1}\equiv 0$. 
So in this case, we have no more the canonically defined subbundles $E^+_2$ and 
$E^-_2$ as before (see Subsection $5.1$). 
Here the underlying geometric structure is 
$$g_4 :=(X, E^+, E^-, \omega).$$
Denote by $G'$ the isometry group of $g_4$. Fix $x\in \widetilde M$ and denote by $H'$ the isotropy 
subgroup of $x$. Then we have $\widetilde M\cong G'/H'.$

To simplify the notations, we identify $T_xM$ with $T_x\widetilde M$. Take a dual basis 
of ${E}^+_x\oplus E^-_x$ with respect 
to ${\omega}_x\mid_{{E}^+_x\oplus {E}^-_x}$, $(y^+, 
l^+_1, l^-_1, y^-)$, such that $l^\pm_1\in {E}^\pm_1$ and $d\lambda(y^+, y^-) =1.$ 
Denote by $\varphi$ the section of $End (E^+\oplus E^-)$, such that 
$$d\lambda(\cdot, \  \cdot ) = \omega (\varphi\  \cdot, \  \cdot).$$ 
Since $d\lambda\wedge\omega =0$, then Tr$\varphi =0.$ So $\exists\  B\not =0,$ such that  
$$\varphi_x =\left (\begin{array}{cccc}
0 & 0 & 0 & 0\\
B & 0 & 0 & 0\\
0 & 0 & 0 & B\\
0 & 0 & 0 & 0 
\end{array}\right ).$$
For $\forall \  h\in H'_e$, $D_xh$ preserves $d\lambda_x$. 
So in the basis above, the matrix of $D_xh$ must be of the following form
$$D_xh =\left ( \begin{array}{cccc}
c & 0 & 0 & 0\\
d & c & 0 & 0\\
0 & 0 & \frac{1}{c} & -\frac{d}{c^2}\\
0 & 0 & 0 & \frac{1}{c}
\end{array}\right )$$
So $H'_e$ is Abelian and isomorphic to $0$, $\mathbb{R}$, 
or $\mathbb{R}^2$. Then as in Subsection $4.1$, $G'_e$ is seen 
to be simply connected.

Denote by $\mathfrak{g}'$ and $\mathfrak{h}'$ the Lie algebras of $G'$ and $H'$. 
Then we get a similar identification of $\mathfrak{g}'$ and $T_xM\oplus\mathfrak{h}'$ 
as in Subsection $5.1.$ In particular, if 
$L_0 :=P^+- P^-$, then $L_0\in \mathfrak{h}'$. We deduce that $\mathfrak{h}' \cong \mathbb{R}$ or 
$\mathbb{R}^2$.

Lemma $4.1.2.$ is still valid here. But the proofs of Lemmas $4.1.3.$ and $5.2.1.$ do not 
pass through in the current case.\\\\
{\bf 6.2. dim$\mathfrak{h}'$=1.}

In this subsection, we suppose that dim$\mathfrak{h}'=1.$ 
Then $\mathfrak{g}' = \mathbb{R}L_0\oplus T_xM.$ By the Jacobi identities of 
$\mathfrak{g}'$, we get the following relations with respect to the dual basis in the 
previous subsection,   
$$[L_0, y^\pm] =\pm y^\pm,\   [L_0, l^\pm_1] = \pm l^\pm_1,$$
$$[X_x, y^\pm] =\mp y^\pm,\   [X_x, l^\pm_1] =\mp l^\pm_1,$$
$$[y^+, y^-] = -X_x - L_0.$$
The brackets, not appeared above, vanish. 

As in Subseciton $5.3$, we define $\alpha :=X_x + L_0$ and 
$$\mathfrak{g} :=\mathbb{R}\alpha
\oplus\mathbb{R}y^+\oplus\mathbb{R}y^-\oplus \mathbb{R}l^+_1\oplus \mathbb{R}l^-_1.$$ 
Thus $\mathfrak{g}'\cong \mathfrak{g}\rtimes \mathbb{R}L_0.$ Denote by $G_e$ 
the connected Lie subgroup of $G'_e$ with Lie algebra $\mathfrak{g}$. Then $G_e$ is naturally 
identified to $\widetilde M$ under $\psi$ (see Subsection $5.3.$) and the corresponding geometric structure
on $G_e$ is given by 
$$g_5 := (\alpha, \mathbb{R}l^+_1\oplus\mathbb{R}y^+, \mathbb{R}l^-_1\oplus\mathbb{R}y^-, 
\psi^\ast\omega ).$$
In addition, $(y^+, l^+_1, l^-_1, y^-)$ is dual with respect to $\psi^\ast\omega.$ 
 
For $\forall\  c,\  d\in \mathbb{R}$, there is a unique Lie algebra automorphism of $\mathfrak{g}$, $\rho^c_d$, 
such that 
$$\rho^c_d(y^\pm) = exp^{\pm c} (y^\pm\pm d\cdot l^\pm_1), $$
$$\rho^c_1( l^\pm_1) =exp^{\pm c}\cdot l^\pm_1, \  \rho^c_d(\alpha) =\alpha.$$
Their corresponding isomorphisms of 
$G_e$ forms a Lie group isomorphic to $\mathbb{R}^2$. Then we 
observe easily that its action on $G_e$ preserves $g_5$ and fixes $e$. So 
dim$( H'_e)\geq 2$, which is a contradiction.\\\\
{\bf 6.3. dim$\mathfrak{h}' $=2.}

In this subsection, we suppose that dim$\mathfrak{h}' =2.$ 

$\ $

{\bf Lemma 6.2.1.} {\it $\exists\  c< 2$, such that 
$\Omega^+ = c \cdot d\lambda.$}

{\it Proof.} Let $\varphi$ and $\psi$ be the same endomorphisms as in 
Subsections $4.1.$ and $6.1.$ Take $l^+_1\in (E^+_1)_x.$ Thus $\varphi_x( l^+_1) =0$. 

Since $\Omega^+\wedge\Omega^+=0$, then det$\psi_x =0.$ 
So if $\psi_x(l^+_1)\not =0$, then there exists $y^+\not =0$, such that 
$\psi_x(y^+) =0.$ Extend $y^+$ and $l^+_1$ to a dual basis, 
$(y^+, l^+_1, z^-, y^-)$. Then in this basis, we get 
$$\varphi_x =\left (\begin{array}{cccc}
a & 0 & 0 & 0\\
b & 0 & 0 & 0\\
0 & 0 & a & b \\
0 & 0 & 0 & 0
\end{array}\right ),\   
\psi_x =\left (\begin{array}{cccc}
0 & A & 0 & 0 \\
0 & B & 0 & 0 \\
0 & 0 & 0 & 0 \\
0 & 0 & A & B
\end{array}\right ). $$
Since $\Omega^+\wedge\omega =0$ and $d\lambda\wedge\omega =0$, then Tr$\varphi$ = Tr$\psi$ =$0$, 
i.e. $a = B =0$. Since $d\lambda\wedge \Omega^+ =0$, then $b\cdot A =0.$ So we get 
$A =0,$ i.e. $\psi_x(l^+_1) =0$, which is a contradiction. 

We deduce that $\psi_x(l^+_1) =0$. Extend $l^+_1$ to a dual basis $(y^+_1, l^+_1, y^-_1, z^-_1)$.  
Thus in this basis, $\varphi_x$ and $\psi_x$ are proportional. Then by homogeneity, we deduce 
the existence of $c$, such that $\Omega^+=c\cdot d\lambda.$

Denote by $\lambda$ the {\it canonical $1$-form} of $\phi_t$ and 
by $J$ be the section of $End(TM)$, such that 
$$J(X) = 0,\  J(u^\pm) =\pm u^\pm.$$
We introduce another $\phi_t$-invariant connection 
$$\bar\nabla :=\nabla -\frac{c}{2}\lambda\otimes J.$$ 
Thus 
$$\bar\nabla_X Y^\pm =[X, Y^\pm] \pm(1-\frac{c}{2}) Y^\pm.$$
Denote by $\bar\Omega^+$ the curvature form of the induced connection $\bar\nabla^+$ 
of $\bar\nabla$ on $\wedge^2 E^+$. Then from the definition of $\bar\nabla$, we easily get 
$$\bar\Omega^+ \equiv 0.$$
Fix a nowhere-vanishing section $\omega^+$ of $\wedge^2 E^+$. Then with respect to 
$\omega^+$, the connection form of $\bar\nabla^+$ is given by 
$$\bar\nabla \omega^+= \bar\beta^+(\cdot)\omega^+.$$
So we have $d \bar\beta^+= \bar\Omega^+ =0.$

Suppose that $c\geq 2$. Then $1 -\frac{c}{2}\leq 0.$ Define $\alpha_t :=\frac{1}{t}\int_{t}^{0}\phi_s^\ast\bar\beta^+ ds.$ 
By the arguments in Subsection $4.4.2.$ of {\bf [BFL2]}, if $t<<0$, then we have  
$$\alpha_t(X) >0.$$
Thus fix $t<<0$ and denote this $\alpha_t$ by $\alpha$. 
Since $\bar \beta^+$ is closed, then so 
is $\alpha$. Define $Y :=\frac{X}{\alpha (X)}$. By Lemma $2.1.2$, the flow 
of $Y$, $\phi^Y_t$, is also a {\it geometric} Anosov flow with smooth distributions. 

Since $\lambda$ is $\phi_t$-invariant, then 
$$ 0 =\mathcal{L}_X\lambda = i_X d\lambda = \alpha (X) (i_Y d\lambda).$$
So 
$$\mathcal{L}_Y d\lambda =i_Y d(d\lambda) +d i_Y d\lambda =0,$$
i.e. $d\lambda$ is $\phi^Y_t$-invariant. Since $\alpha$ is easily seen 
to be the {\it canonical 1-form} of $\phi^Y_t$ and $d\alpha =0$, then 
rank$(\phi^Y_t) =0$. Thus by Subsection $3.1$, $\phi^Y_t$ is 
finitely covered by the suspension of a hyperbolic automorphism of $\mathbb{T}^4$, 
which is given 
by a hyperbolic matrix in $GL(4, \mathbb{Z})$. Then by a direct calculation, 
using the Jordan form of this matrix, $\lambda$ is seen to be closed (see 
{\bf [Fa]} for the details). We deduce that 
rank$(\phi_t) =0$, which is a contradiction.  
$\square$

$\ $

We can see that Lemma $3.2.3.$ is also true for $\bar\nabla$ defined in the previous lemma. 
So in particular, we get  
$$\bar\nabla \bar R =0, \  \bar\nabla \bar T =0,\  \bar T(E_a, E_b)\subseteq E_{a+b},$$
and if $a+b\not =0,$ $\bar R( E_a, E_b) =0,$ where $\bar T$ and $\bar R$ are the torsion and 
curvature tensors of $\bar\nabla.$ Since the $\bar\nabla$-geodesics tangent to 
$E^+$ or $E^-$ are also complete, then by Lemma $A$ in the appendix, $\bar\nabla$ is complete. 
Thus as in Subsection $4.1$, we get the following identification via $\bar\nabla$
$$\mathfrak{g}'\stackrel{\sim}{\longmapsto}T_x\widetilde M\oplus\mathfrak{h}'$$
$$ u\to (Y^u(x),\  ({\widetilde{\bar\nabla}}_{Y^u}-\mathcal{L}_{Y^u})\mid_x).$$

Since $\bar\Omega^+\equiv 0$, then we can define a character $\bar\chi$ of $\mathfrak{g}'$ as in Subsection 
$4.1.$ Thus $\bar\chi^{-1}(0)$ is an ideal of $\mathfrak{g}'$, denoted again by $\mathfrak{g}$. By the same type of 
arguments as before, we easily get
$$\mathfrak{g}\cong \mathbb{R}^3\rtimes \mathfrak{sl}(2,\mathbb{R}).$$

Denote by $G_e$ the connected Lie subgroup of $G'_e$ with Lie algebra $\mathfrak{g}$ and define 
$H_e := H'_e\cap G_e.$ Since $G'_e$ is simply connected, then so is $G_e$ (see Subsection $4.1.$ 
of {\bf [BFL1]}). Thus by some direct calculations, $G_e$ and $H_e$ can be realized as following
$$G_e\cong \mathbb{R}^3\rtimes \widetilde{SO_0(1, 2)},$$
where ${SO_0(1, 2)}$ is the identity component of 
the isometry group of the quadratic form : $-dx^2 + dy^2 + dz^2.$ The semi-direct product 
is given by the composition 
of the projection of $\widetilde{SO_0(1, 2)}$ onto $SO_0(1, 2)$ and the linear 
action of $SO_0(1, 2)$ on $\mathbb{R}^3$. 
Let $((0, 0, 1), 0)\in \mathbb{R}^3\rtimes
\mathfrak{so}(1, 2).$ Then $H_e$ is just the $1$-parameter subgroup 
generated by this vector, denoted also by $\mathbb{R}$. 

Since $\bar\Omega^+\equiv 0$, then the same argument as in 
Subsection $4.2.$ of {\bf [BFL1]} works in our case, if 
we replace the metric entropy there by $2(1-\frac{c}{2}).$ 
Thus the general argument of Section $5.$ 
of {\bf [BFL1]} gives a discrete subgroup of 
$G_e$, acting freely, properly and cocompactly on 
$G_e/H_e$. Now we finish the proof 
by showing

$\ $

{\bf Lemma 6.2.2.} {\it $\mathbb{R}^3\rtimes\widetilde{SO_0(1, 2)}$
admits no discrete subgroup, which acts properly, freely and cocompactly on 
$(\mathbb{R}^3\rtimes \widetilde{SO_0(1, 2)})/\mathbb{R}.$}

{\it Proof.} Recall that $\mathbb{R}$ denotes $H_e$ and $G_e$ denotes $\mathbb{R}^3\rtimes
\widetilde{SO_0(1, 2)}$. Suppose the existence of a subgroup $\Gamma$ satisfying 
the conditions in the lemma. Denote by $\bar\Gamma$ the Zariski closure of 
$\Gamma$ in $G_e$ (Here the Zariski topology of $G_e$ means the lifted topology of the Zariski 
topology of $\mathbb{R}^3\rtimes SO_0(1, 2)$ by the canonical projection). 

If $\Gamma$ is solvable, then $\bar\Gamma$ is also solvable. Then by 
{\bf [Ra]}, there exists a connected closed subgroup $H$
of $\bar\Gamma$, such that $\Gamma\subseteq H$ and $H/\Gamma$ is 
compact. Let cd$(\cdot)$ denote the cohomological dimension of a group. Since 
$\Gamma$ acts cocompactly on $G_e/\mathbb{R}$, then cd$(\Gamma) =5$. We 
deduce that cd$(H) =5$. So $H$ is a closed solvable subgroup of $G_e$ of dimension 
$5$. Then the single possibility is $\mathbb{R}^3\rtimes AN$ ( where $KAN$
is the Iwasawa decomposition of $\widetilde{SO_0(1, 2)}$). But 
$\mathbb{R}^3\rtimes AN$ is not unimodular. So 
it has no cocompact lattice. We deduce that $\Gamma$ is not solvable. Then $\bar\Gamma$
must contain $\widetilde{SO_0(1, 2)}$. 

Since $\Gamma$ acts cocompactly on
$G_e/\mathbb{R}$,  then $\widetilde{SO_0(1, 2)}\subsetneq \bar\Gamma$. We deduce that 
$\bar\Gamma \cap \mathbb{R}^3\not =0$. Since the representation 
of $\widetilde{SO_0(1, 2)}$ on $\mathbb{R}^3$ is irreducible, then $\bar\Gamma$ must be $G_e$, 
i.e. $\Gamma$ is Zariski-dense in $G_e$. Let $\Delta$ be the projection of $\Gamma$
into $\widetilde{SO_0(1, 2)}$, then by [$\mathbf{Ra}$], $\Delta$ is discrete in
$\widetilde{SO_0(1, 2)}$. We deduce that $\Gamma\cap\mathbb{R}^3\not =0$. 
Since the semi-direct product is given by an irreducible representation, $\Gamma
\cap\mathbb{R}^3$ is in fact cocompact in $\mathbb{R}^3$.

Since $\Gamma$ acts properly on $G_e/\mathbb{R}$, then $\Gamma\cap\mathbb{R}^3$ acts 
properly on $\mathbb{R}^3/\mathbb{R}$ which is 
a closed subset of $G_e/\mathbb{R}$. We deduce that $\mathbb{R}^3$ 
acts also properly on $\mathbb{R}^3/\mathbb{R}$. But it is absurd.  $\square$\\\\
{\bf Appendix.}

At first, we prove  the following elementry lemma, which is used in the proof of 
Lemma $3.2.4$.

$\ $

{\bf Lemma A.} {\it Let $\nabla$ be a smooth linear connection on a 
connected manifold $M$ of dimension $n$. Let $X_1, \cdots, X_k$ be complete 
fields on $M$ and $E_1, \cdots, E_l$ be smooth distributions on $M$, such that \\
$(1)$. $\nabla X_i =0, \forall \  1\leq i\leq k,$ $\nabla E_j\subseteq E_j, 
\forall \  1\leq j\leq l, $\\
$(2)$. $TM = \mathbb{R}X_1\oplus\cdots\oplus \mathbb{R}X_k
\oplus E_1\oplus\cdots\oplus E_l ,$\\
$(3)$. $\nabla R =0, $ $\nabla T =0,$\\
$(4)$. $\forall \  1\leq j\leq l$, the geodesics of $\nabla$, tangent to $E_j$, are 
defined on $\mathbb{R}$, \\
then $\nabla$ is complete.}

{\it Proof.} For the terminology below, our reference is {\bf [K-No]}, vol. I. 
For $\forall \  1\leq i\leq k$, since $X_i$ is complete and parallel, then 
any geodesic tangent to $\mathbb{R}X_i$ is defined on $\mathbb{R}$. 
So without any loss of generality, we suppose that $k =0.$

Let $\mathcal{F}(M)$ be the frame bundle of $M$ and $\pi$ the projection of 
$\mathcal{F}(M)$
onto $M$. The linear connection $\nabla$ gives a horizontal distribution 
$\mathcal{H}$ on $\mathcal{F}M$ and $\mathcal{F}M$ is foliated by holonomy 
subbundles. $\mathcal{H}$ is tangent to each holonomy subbundle, then so 
is any standard horizontal field. For $\forall\  u \in \mathcal{F}M$, denote by $P(u)$ the 
holonomy subbundle containing $u$. The induced fields on $P(u)$ of the standard
horizontal fields of $\mathcal{F}M$ are called also standard horizontal. By 
{\bf [K-No]}, $\nabla$ is complete, iff for $\forall\  x\in M, \exists\  u\in \pi^{-1}
(x)$, such that the standard horizontal fields of $P(u)$ are all complete.

Take $x\in M$ and $u\in \pi^{-1}(x)$, such that
$$ u =( v^1_1,\cdots, v^1_{i_1}, \cdots, v^l_1, \cdots, v^l_{i_l} ),$$
where $\{ v^j_1, \cdots , v^j_{i_j} \}$ is a basis of $E_j(x)$, $\forall \  1\leq j\leq l.$ 
For $\forall\   \xi\in \mathbb{R}^n$, the standard horizontal field on $P(u)$ corresponding to $\xi$ 
is denoted by $B^u(\xi)$ and the canonical basis of $\mathbb{R}^n$ is denoted by 
$(e_1,\cdots, e_n )$. Take $v\in P(u)$. Because of assumption $(1)$, 
$v$ has the same form as $u$. Then for $\forall$ $ 1\leq m\leq n$, the integral curve 
of $B^u( e_m)$, begining at $v$, is just the horizontal lift, beginning at $v$, of the 
geodesic tangent to $Pr_m( v )$. By assumption $(4)$, such a geodesic 
is defined on $\mathbb{R}$. We deduce that $B^u( e_m)$ is complete. 

Fix a basis of the holonomy algebra of $\nabla$ and denote the corresponding 
vertical fields of $P(u)$ by $\{ V_1, \cdots , V_s \}$. By assumption $(3)$, 
the fields $$\{ V_1, \cdots , V_s, B^u( e_1), \cdots , B^u( e_n ) \}$$
 generate a 
Lie algebra. Since these fields are all complete, then this Lie algebra must be 
induced by the smooth action on $P(u)$ of a simply connected Lie 
group. Thus for $\forall\  \xi\in \mathbb{R}^n$, the field 
$B^u(\xi )$ $( = \sum_{1\leq i\leq n}\xi_i B^u( e_i ) )$ is complete. We 
deduce that $\nabla$ is complete.  $\square$

$\ $

The following lemma is used in the proof of Lemma $2.1.2.$

$\ $

{\bf Lemma B.} {\it Let $\phi_t$ be an Anosov flow with $C^\infty$ distributions on 
a closed manifold $M$. If $f$ is a smooth positive function on $M$ and the flow 
of $fX$($X$ is the generator of $\phi_t$) has also $C^\infty$ 
distributions, then there exists a $C^\infty$ $1$-form $\alpha$ on $M$, such that 
$\mathcal{L}_X d\alpha =0$ and $f =\frac{1}{\alpha(X)}$. Conversely, if $\alpha$ is a 
$C^\infty$ $1$-form on $M$, such that $\mathcal{L}_X d\alpha =0$ and $\alpha(X) >0$, 
then the flow of $\frac{X}{\alpha(X)}$ has also $C^\infty$ distributions.} 

{\it Proof}. Recall at first that a $C^\infty$ time change of an Anosov flow is also Anosov. 
Let $fX$ be a time change of $\phi_t$ with smooth distributions. 
Denote by $\phi^{fX}_t$ the flow of $fX$ 
and by $\lambda_1$ its {\it canonical 1-form}. 
Then $\lambda_1(fX) =1$, i.e. $f = \frac{1}{\lambda_1(X)}$. 
Since $\lambda_1$ is $\phi^{fX}_t$-invariant, then $i_{fX}d\lambda_1 =0$. So 
$i_X d\lambda_1 =0$. We deduce that $\mathcal{L}_X d\lambda_1 =0.$

If $hX$ is a smooth time change of $\phi_t$, then its strong stable 
distribution is given by (see Lemma $1.2.$ of {\bf [LMM]})
$$ E_{hX}^- =\{ Y^- + \beta (Y^-) X \mid Y^- \in E^-_X\},$$
where $E_{hX}^-$ denotes the strong stable distribution of $hX$ and 
$\beta$ is a $C^0$ section of $(E_X^-)^\ast$, such that 
$$\mathcal{L}_{X}(h^{-1}\beta) = h^{-2} dh. \eqno (\ast)$$
Denote by $\lambda$ the {\it canonical 1-form} of $\phi_t$. If $\alpha$ is a smooth 
$1$-form on $M$, such that $\mathcal{L}_X d\alpha =0$ and 
$\alpha(X) >0$, then by a simple calculation, $-\frac{\alpha -\lambda}{\alpha(X)}$ satisfies 
the previous equation $(\ast)$ about $\beta$ with $h :=\frac{1}{\alpha(X)}$. So 
$E^-_{\frac{X}{\alpha(X)}}$ is smooth. Similarly, $E^+_{\frac{X}{\alpha(X)}}$ 
is also smooth.  $\square$      

$\ $

{\bf Acknowledgements.} The author would like to thank his thesis advisors, P. Foulon 
and P. Pansu, for the discussions and help. He would like also to thank Y. Benoist and 
F. Labourie for their help.

$\ $
     
{\bf References}

$\ $

{\small 
{\bf [Am]} A. M. Amores, Vector fields of a finite type G-structure, {\it J. Differential Geometry} 14 (1979) 1-6.

{\bf [An]} V. D. Anosov, Geodesic flows on closed Riemannian manifolds with negative curvature, 
{\it Proc. Inst. Steklov} 90 (1967) 1-235.

{\bf [Ba]} T. Barbot, Caract\'erisation des flots d'Anosov en dimension 3 par leurs feuilletages faibles, 
{\it Ergod. Th. and Dynam. Sys.} 15 (1995) 247-270.

{\bf [BFL1]} Y. Benoist, P. Foulon and F. Labourie, 
Flots d'Anosov \`a distributions de Liapounov diff\'erentiables. I,  
{\it Ann. Inst. Henri Poincar\'e} 53 (1990) 395-412.

{\bf [BFL2]} Y. Benoist, P. Foulon and F. Labourie, 
Flots d'Anosov \`a distributions stable et instable diff\'erentiables, 
{\it J. Amer. Math. Soc.} 5 (1992) 33-74.

{\bf [BL]} Y. Benoist and F. Labourie, Sur les diff\'eomorphismes
d'Anosov affines \`a feuilletages stable et instable diff\'erentiables, 
{\it Invent. Math.} 111 (1993) 285-308.

{\bf [C-Q]} A. Candel and R. Quiroga-Barranco, Gromov's centralizer theorem, to appear in {\it Geometriae Dedicata.}

{\bf [Fa]} Y. Fang, A remark about hyperbolic infranilautomorphisms, 
{\it C. R. Acad. Sci. Paris, Ser. I 336 No.9} (2003) 769-772.

{\bf [FK]} R. Feres and A. Katok, invariant tensor fields of dynamical systems 
with pinched Lyapunov exponents and rigidity of geodesic flows, {\it 
Ergod. Th. and Dynam. Sys.} 9 (1989) 427-432.

{\bf [Gh]} \'E. Ghys, Flots d'Anosov dont les feuilletages stables sont 
diff\'erentiables, {\it Ann. Scient. \'Ec. Norm. Sup. (4)} 20 (1987) 251-270.

{\bf [HK]} S. Hurder and A. Katok, Differentiability, rigidity and 
Godbillon-Vey classes for Anosov flows, {\it Pub. I.H.\'E.S.} 72 
(1990) 5-61.

{\bf [HaK]} B. Hasselblatt and A. Katok, Introduction to the modern theory 
of dynamical systems, {\it Cambridge University Press,} 1995.

{\bf [K]} M. Kanai, Geodesic flows of negatively curved manifolds with 
smooth stable and instable foliations, {\it Ergod. Th. and Dynam. Sys.} 8 (1988) 215-240.

{\bf [K-No]} S. Kobayashi and K. Nomizu, Foundations of Differential 
Geometry, vol. I, II, {\it Interscience, New York and London,} 1963.

{\bf [Lich]} A. Lichn\'erowicz, Th\'eorie globale des connexions et des groupes 
d'holonomie, {\it Edizioni Cremonese, Roma,} 1962. 

{\bf [LMM]} R. de la llave, J. Marco and R. Moriyon, Canonical perturbation theory of 
Anosov systems and regularity results for Livsic cohomology equation, {\it Ann. Math} 
123(3) (1986) 537-612.

{\bf [Plan]} J. F. Plante, Anosov flows, {\it Amer. J. Math.} 94 (1972) 729-754.

{\bf [Ra]} M. S. Raghunathan, Discrete subgroups of Lie groups, 
{\it Springer, Berlin, Heidelberg, New York,} 1972.

{\bf [To]} P. Tomter, Anosov flows on Infra-homogeneous Spaces, {\it Proc. Symp. in Pure Math, 
Vol. XIV, Global Analysis,} (1970) 299-327.}



\end{document}